\documentclass{article}[12pt]
\usepackage{color,times,amsmath,amsfonts,latexsym,epsfig,epsf,colordvi}
\usepackage{url,hyperref}
\usepackage[english]{babel}

\newcommand{\CC}{\mathbb {C}}
\newcommand{\RR}{\mathbb {R}}

\newcommand{\al}{\alpha }

\newcommand{\bp}{\begin{pmat}}
\newcommand{\ep}{\end{pmat}}

\author{  Frank Uhlig \thanks{Department of Mathematics and Statistics, Auburn 
University, Auburn, AL 36849-5310 \ (uhligfd@auburn.edu)}}

 \title{\vspace*{-12mm} Constructing   the Field of Values of  Decomposable and General Matrices Using the ZNN Based Path Following Method\\}
\begin{document}
\date{~}
\thispagestyle{empty}
\maketitle

\thispagestyle{empty}

\vspace*{-16mm}
\begin{center} { \bf Abstract  } \\[2mm]
\begin{minipage}{140mm}
This paper  describes and develops a fast and accurate  path following algorithm that computes the field of values boundary curve $\partial F(A)$ for every conceivable complex or real square matrix $A$. It relies on the matrix flow decomposition algorithm that finds a proper block-diagonal flow representation for the associated hermitean matrix flow ${\cal F}_A(t) = \cos(t) H + \sin(t) K$ under unitary similarity if that is possible. Here ${\cal F}_A(t)$  is the 1-parameter-varying linear combination of the real and skew part matrices $H = (A+A^*)/2$ and $K = (A-A^*)/(2i)$ of $A$. For indecomposable matrix flows, ${\cal F}_A(t)$ has just one block and the ZNN based field of values algorithm works with ${\cal F}_A(t)$ directly. For decomposing  flows ${\cal F}_A(t)$, the algorithm decomposes the given  matrix $A$ unitarily  into block-diagonal form $U^*AU =  \text { diag} (A_j)$  with $j > 1$ diagonal blocks $A_j$ whose individual  sizes add up to the  size of  $A$. It  then computes the field of values boundaries separately  for each diagonal block $A_j$ using  the  path following ZNN eigenvalue method.  The convex hull of all sub-fields of values boundary points $\partial F(A_j)$ finally determines the field of values boundary curve correctly for  decomposing matrices $A$. The algorithm removes standard restrictions for path following FoV methods that generally cannot deal with decomposing matrices $A$ due to possible eigencurve crossings of ${\cal F}_A(t)$. Tests and numerical comparisons are included. Our ZNN based method is coded for sequential and parallel computations and both versions run very accurately and   fast when compared with Johnson's  Francis QR eigenvalue  and Bendixson rectangle based method  and compute global eigenanalyses of ${\cal F}_A(t_k)$ for  large discrete sets of angles $t_k  \in {[} 0,2\pi{]}$ more slowly.
\end{minipage}\\[-1mm]
\end{center}  
\thispagestyle{empty}

\hspace*{17.5mm}  decomposable matrix, numerical algorithm,  block-diagonal matrix, unitary similarity\\[-3mm]


\pagestyle{myheadings}
\thispagestyle{plain}
\markboth{Frank Uhlig}{Decomposing and General Matrix FoVs }

\section{Field of Values Computations for  Real or Complex Square Matrices; Preliminaries}\vspace*{-1mm}

This paper develops an algorithm to construct the {\em Field of Values} (FoV) of a  real or complex square matrix $A$ via any chosen path following method, even when $A$ is unitarily block-diagonalizable into $j > 1$ diagonal blocks $A_j$. Constructing the FoV  of properly  decomposable matrices $A$ was previously only possible by using {\em global matrix eigenvalue finders} such as Francis' QR algorithm that compute all eigendata of a matrix. Here we solve this problem in three steps by \\[0.5mm]
  (a) checking for unitary diagonalizability of a given matrix $A$ via ${\cal F}_A(t)$, see \cite {FUDecomp},  and unitarily block-diagonalizing  decomposable matrices into $j >1$ diagonal blocks, then\\[0.5mm]
   (b) using a ZNN based path following method to compute discrete FoV boundary eigendata  for each  diagonal block $A_j$ separately, and finally by \\[0.5mm]
   (c) implementing a convex hull algorithm to plot the FoV boundary curve of $A$ itself.\\[2mm]
The  {\em field of values} $F(A)$, also known as the {\em numerical range}  of a real or complex square matrix $A$ is defined as
$$  F(A)\ = \ \{ x^*Ax \mid x \in \CC^n, \|x\|_2 = 1\}  \subset \CC \ .$$
Its genesis, original purpose  or inventors are unknown. In 1918/1919 Hausdorff and Toeplitz showed independently that $F(A)$ is  convex for any square matrix $A$. Other well known properties of the field of values $F(A)$ are:\\
The field of values $F(A)$ is compact for any $A$; the field of values  of a 2 by 2 matrix is an ellipse; $F(A)$ contains all eigenvalues of $A$; $F(c A) = cF(A)$ for all scalars $c \in \CC$; $F(U^*AU) = F(A)$ for any unitary matrix $U$; for normal matrices $A$, $F(A)$  is the convex hull of $A$'s eigenvalues; for block-diagonal matrices 
$$A=\bp A_1 & O\\O & A_2\ep \ , \ \ F(A) = \text{ convhull }  ( F(A_1), F(A_2)) .
$$
  The matrix field of values is treated in many textbooks, such as in \cite{HJ94}. The maximal and minimal distant points of the boundary $\partial F(A)$ of the field of values from the origin $0 \in \CC$ play an important role in matrix norm estimations and in control theory, respectively. The former is called the {\em numerical radius} of the matrix and the latter is the {\em Crawford number} of $A$. Generalized numerical ranges are important notions in the study of quantum computations and elsewhere.\\[2mm]
 Building on Bendixson rectangles \cite{Be} from 1902 that contain the field of values of a matrix, Johnson  \cite{J78} in 1978 established an eigendata based method to find  discrete boundary curve points on $\partial F(A)$ of the field of values  for a matrix $A$ via repeated hermitean eigendata computations. Johnson's method involves the  matrices\\[-2mm]
  $${H=  (A+A^*)/2 = H^*} \ \ \text{ and } \ \
 {K=  (A - A^*)/(2i) =  K^*}  \in \CC_{n,n} $$
   with $A = H + iK$ and the  associated 1-parameter-varying   hermitean matrix flow\\ 
 \begin{equation}
 \label{Aoft} {\cal F}_A(t) = \cos(t)  H +   \sin(t)   K  = ({\cal F}_A(t))^*  \ \ \ \text{ for  angles } \ \ \ 0 \leq t \leq 2\pi \ .
\end{equation}  
More specifically, the normalized  eigenvectors $x(t)$ and $y(t)$ for the largest and smallest eigenvalues of each ${\cal F}_A(t)$ determine two $\partial F(A)$ { boundary points} via the quadratic form evaluations  $x(t)^*Ax(t) \text{ and } y(t)^*Ay(t) \in \CC$  and two $\partial F(A)$ tangents. One standard approach to approximate  the boundary curve $\partial F(A)$ of a given $n$ by $n$  matrix $A$ with $n$ up to around $n=10,000$  is to use the Francis implicit multi-shift QR algorithm  for many discrete parameter values $0 \leq t_i \leq 2\pi$ and find the extreme eigendata of each ${\cal F}_A(t_i)$. Then evaluate  the quadratic forms $x(t_i)^*Ax(t_i) \text{ and } y(t_i)^*Ay(t_i)$ for the respective 'extreme' eigenvectors $x(t_i)$ and $y(t_i)$ of ${\cal F}_A(t_i)$ in order to  find accurate  interpolation points on the FoV boundary curve   $\partial F(A)$. For huge dimensions $n \gg 10,000$ other ways can be used to find the extreme eigenvalues and associated eigenvectors iteratively such as Matlab's {\tt eigs}  m-file.\\[1mm]
For relatively small dimensioned matrices $A$, the global hermitean QR eigenvalue approach computes complete eigendata  sets for  each ${\cal F}_A(t_i)$ reliably.  {Note that for each angle $t_i$ in Johnson's method we need to compute all eigenvalues of ${\cal F}_A(t_i)$ at  large computational expense because the eigencurves of ${\cal F}_A(t)$ may suffer lead changes or eigencurve crossings. Only the locally  extreme ones help us to  plot $\partial F(A)$ boundary points. Therefore non-global but local path following methods have become attractive recently.} Francis' implicit multi-shift QR  and other complete matrix eigensolvers offer a foolproof way to solve the FoV problem for all matrices $A$ of modest size, decomposable or not. For large dimensions $n$, with $n$ near $10,000$, global eigen methods become rather  expensive when compared with  simpler and much faster path following methods.\\
 But until now, path following methods alone could not handle the FoV problem for unitarily decomposable matrices,  see e.g, \cite[Algorithm 6.1 and p. 1733 - 1743]{LM18} for a detailed analysis. How to determine whether a given complex matrix $A$ can be block-decomposed by a fixed unitary matrix similarity  has recently been solved constructively in \cite{FUDecomp}. Attempts at unitary block-decompositions of static entry matrices and general matrix flows go back at least 90 years to the beginnings of quantum physics and quantum chemistry when von Neumann and Wigner \cite{NW29}  established two (in-)decomposability criteria for 1-parameter hermitean matrix flows from their eigencurve behavior. Namely, if for some $0 \leq t < 2\pi$ two eigencurves for ${\cal F}_A(t)$  cross each other, then their respective eigendata is associated with  separate diagonal blocks of the hermitean matrix flow  ${\cal F}_A(t)$ -- and thus by \cite{FUDecomp} also of $A$. If on the other hand two eigencurves seem to be attracted to each other and almost touch but don't actually cross but veer off like the two branches of a hyperbola do, then the respective eigencurves belong to the same unitarily indecomposable block of $A$, see \cite{NW29} and \cite{FUDecomp} for further details.\\
 {The previously unsolvable {\em block diagonalization problem} for real or complex static matrices $A$ under  unitary similarity has been solved recently and constructively in \cite{FUDecomp} by looking at two different 'field of values' function matrices ${\cal F}_A(t_1)$ and ${\cal F}_A(t_2)$. Diagonalize the hermitean matrix ${\cal F}_A(t_1)$ via QR and perform a matrix similarity with the ONB eigenvector matrix $U_1$ for  ${\cal F}_A(t_1)$ on ${\cal F}_A(t_2)$. If after zeroing out all entries in $W = U_1^*{\cal F}_A(t_2) U_1$ of magnitudes below $10^{-14} ||A||$ there are multiple zero entries in $W$, sort the eigenvectors in $U_1$ according to the 0-1 pattern of $W$ and  $A$ is block diagonalized by the unitary similarity with the thus amended $U_1$, see \cite{FUDecomp}. } \\[2mm]
Some path following methods for 1-parameter or time-varying matrix problems  {- unlike ours -} track the solution of a derived differential equation by using numerical integrators. To do so for the matrix FoV problem, Loisel and Maxwell \cite[sections 4 and 6]{LM18}  for example,  differentiate the Bendixson/Johnson  eigenvalue equation 
$$ {\cal F}_A(t) u(t) = \lambda(t)u(t) \ $$
with respect to $t$.  
When additionally  requiring  that the eigenvectors $u(t)$ be normalized, they obtain a 1-parameter matrix and vector differential equation in $\CC^{n+1}$ for the matrix FoV problem. Their  ODE path following method uses the eigendata of ${\cal F}_A(0)$ as the initial value and then proceeds with the Dormand--Prince RK5(4)7M numerical integrator \cite{DP80} of fifth order accuracy $O(h^5)$. This  adaptive method can be implemented inside Matlab as {\tt ode45}. {\tt ode45} is quite versatile, accurate  and well suited  for stiff differential equations. The computed $\partial F(A)$ points are then smoothed through  fourth degree Hermite interpolation in \cite{LM18}. Unfortunately the codes and error threshold settings etc used by Loisel and Maxwell in \cite{LM18} are no longer available. For unitarily indecomposable matrices $A$, the path following results  in \cite{LM18} are as accurate as those from the Francis QR matrix eigenvalue algorithm { and more economical} than a complete global Johnson QR {\tt eig} implementation.\\
      The accuracy and speed of computing  indecomposable matrix FoVs was subsequently bested  again in \cite{FUZNNFoV}  by using {\em Zhang Neural Network}  (ZNN) based eigen-computations as the  path following method.\\ {\em Zhang Neural Networks}   are specially designed  to solve time-varying matrix problems {via one linear equations solve and a simple vector recursion per time step}. They proceed in a totally different  way than classical integrators and start from a time-varying matrix model. {Then the model's global error function $E(t)$ is stipulated to decay   exponentially fast in  $t$. Subsequently the postulated error decay differential equation $\dot E(t) = -\eta E(t)$ with $\eta > 0$ is solved for the derivatives of all unknowns $x(..)$ algebraically.} {Then a look-ahead and convergent finite difference formula involving $x(t_{k+1})$ and earlier solutions at times $t_k, t_{k-1}, ...$ and $\dot x(t_k)$ only is chosen from a table of such predictive and convergent finite difference formulas according to the desired truncation error order. And as before, this difference formula is solved for $\dot x(t_k)$ and the two expressions for $\dot x(t_k)$ are equated. This results in a derivative free formula for the future time solution $x(t_{k+1})$ in terms of a linear equations solve and a finite difference formula at each time step. A more detailed derivation and a dozen standard time-varying matrix flow problems are described in the survey article on ZNN methods in \cite{FUZNNsurvey}.}\\[1mm]  
 {Thus the given time-varying matrix problem is solved at each discrete time step $t_k$, giving us the predicted solution $x(t_{k+1})$ in the future  for the matrix flow problem  by using one linear equations solve and a convergent look-ahead finite difference formula evaluation. The difference formula relates the solution $x(t)$ at $t = t_{k+1}$ to earlier known systems data. Practical and  theoretical constructions  of discretized ZNN methods for various time-varying matrix problems are the subject of \cite{FUZNNsurvey}.}\\[2mm]
Today  general path following methods  have become suitable for solving general time-varying matrix problems  due to our ability to decipher dense matrices and matrix flows that are unitarily block-diagonalizable,  both theoretically and computation-wise as detailed in \cite{FUDecomp}. 
 The decomposability check  is elementary, accurate and fast. It works universally for any real or complex square matrix $A$. To check on the unitary block-diagonalizability of a dense matrix $A \in \CC_{n,n}$  with $n = 250$ for example  takes around 0.05 seconds for both decomposable and indecomposable matrices $A$. The new matrix decomposition algorithm \cite{FUDecomp} {is not affected by} the eigen or Jordan structure of  $A$, nor by any other conceivable property or defect of $A$. It derives solely from the properties of the associated hermitean matrix flow ${\cal F}_A(t)$.  And irrespective of  whether $A$ is found to be  unitarily block-decomposable or not by the algorithm in \cite{FUDecomp}, the field of values boundary curve $\partial F(A)$ data can now be computed by using - for example -  any ODE path follower or the fast ZNN method \cite{FUZNNFoV} on  $A$ alone or for each  diagonal block  $A_j$ of a unitary block-diagonalization of $A$ separately  if $A$ is properly decomposable, respectively. To plot the overall FoV boundary curve  in the decomposable case, we use Matlab's {\tt convhull} algorithm on the totality of all computed diagonal block FoV boundary points. This method is  new and much quicker than using a global eigensolver in either case, i.e.,  independently of whether $A$ unitarily decomposes or not.\\[1mm]
This paper combines our recent knowledge of how to  block-diagonalize both static matrices and matrix flows by unitary similarity  \cite{FUDecomp} (which had been unrealized for almost a century) with the recent introduction  of a path following ODE method \cite{LM18} to compute the field of values of indecomposable matrices.  We develop  and test one algorithm that can now find the field of values boundary for unitarily decomposable --  unavailable before -- and indecomposable matrices quickly using a path-follower that is powered by  Zhang Neural Networks  for the single-parameter-varying matrix FoV problem with high accuracy and speed.\\[2mm] 
Section 2  illustrates the new computational landscape for any unitarily invariant matrix problem and shows that decomposable matrices generally do not lend themselves well in unitarily invariant matrix and  'divide and conquer'  computations.  Section 3  explains the complexity gains of path following methods for decomposing matrices and  {further elaborates some advantages and peculiarities} of the ZNN method \cite{FUZNNsurvey}, followed by  numerical test results for our combined FoV algorithm.  Section 4  looks back in history and  forward.\\[-6mm]

\section{Path Following Methods for  Field of Values Computations; a Warning Illustration}

In this section we describe the ill results of trying speedy path following methods to compute the field of values of a  dense matrix $A$ that is unitarily block-diagonalizable.\\
 For this purpose we construct a dense block-diagonalizable 15 by 15 matrix $A$ beginning with  the Matlab command  {\tt   B = blkdiag(1i*randn(10),-randn(5)-(3-2i)*eye(5));} . $B$ is a 15 by 15 non-normal block-diagonal matrix. Then we  create a dense 15 by 15 unitary random entry matrix $Q$ via {\tt [Q,R]= qr(randn(15));} and  form the test matrix  {\tt A = Q'*A*Q;} . Clearly the resulting matrix A is non-normal, dense  and its unitary block-diagonalizability is  hidden from view. Next we  choose any path following method to plot a discrete sample of FoV boundary points from a partial eigenanalysis of ${\cal F}_A(t_i)$ with $ 0 \leq t_i \leq 2\pi$ of $A$ by using quadratic form evaluations with $A$ as described earlier. Which path following method we choose, be it an initial value ordinary differential equations solver such as recently done in \cite{LM18}, or a ZNN based one, see \cite{FUZNNFoV} e.g., or any other method is irrelevant here. Look at the perturbing `FoV  boundary curve' drawing for our dense test matrix $A$  in Figure 1. \\[-6mm]
\begin{center}
\includegraphics[width=95mm]{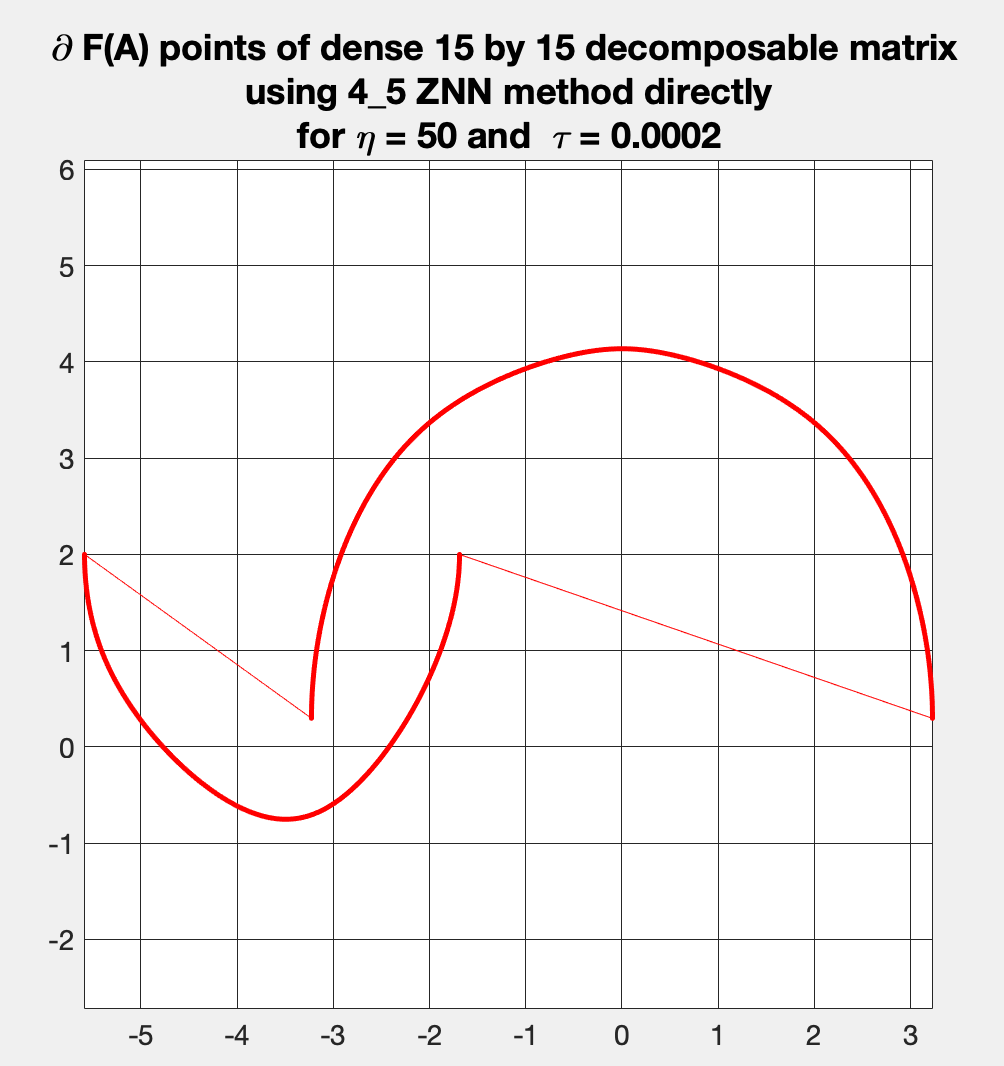} \\[-1mm]
Figure  1: Path following FoV output of a unitarily block-diagonalizable matrix $A$, showing the both extreme eigenvalue induced $\partial F(A)$ points for $0\leq  t \leq \pi$ of ${\cal F}_A(t)$, moving counterclockwise  separately and resulting in two partial and disjoint curves
\end{center}
\vspace*{-0.5mm}
Disjoint FoV mis-representations such as depicted in Figure 1 or similar grotesquely incomplete FoV representations are standard with  path following methods when applied unwittingly to  unitarily block-diagonalizable matrices. By design, these methods follow the initially chosen extreme ${\cal F}_A(t)$ eigenvalue curve faithfully and reflect the problems with computing the FoV boundary naively for decomposable matrices whose extreme eigencurves cross.\\[1mm] 
Next we plot the two individual block FoV boundary curves completely  for our $A$  that are hinted at in Figure 1 after  computing the underlying unitary block decomposition of $A$ via the block-diagonalizing algorithm of \cite{FUDecomp}.\\[-5.5mm]
\begin{center}
\includegraphics[width=95mm]{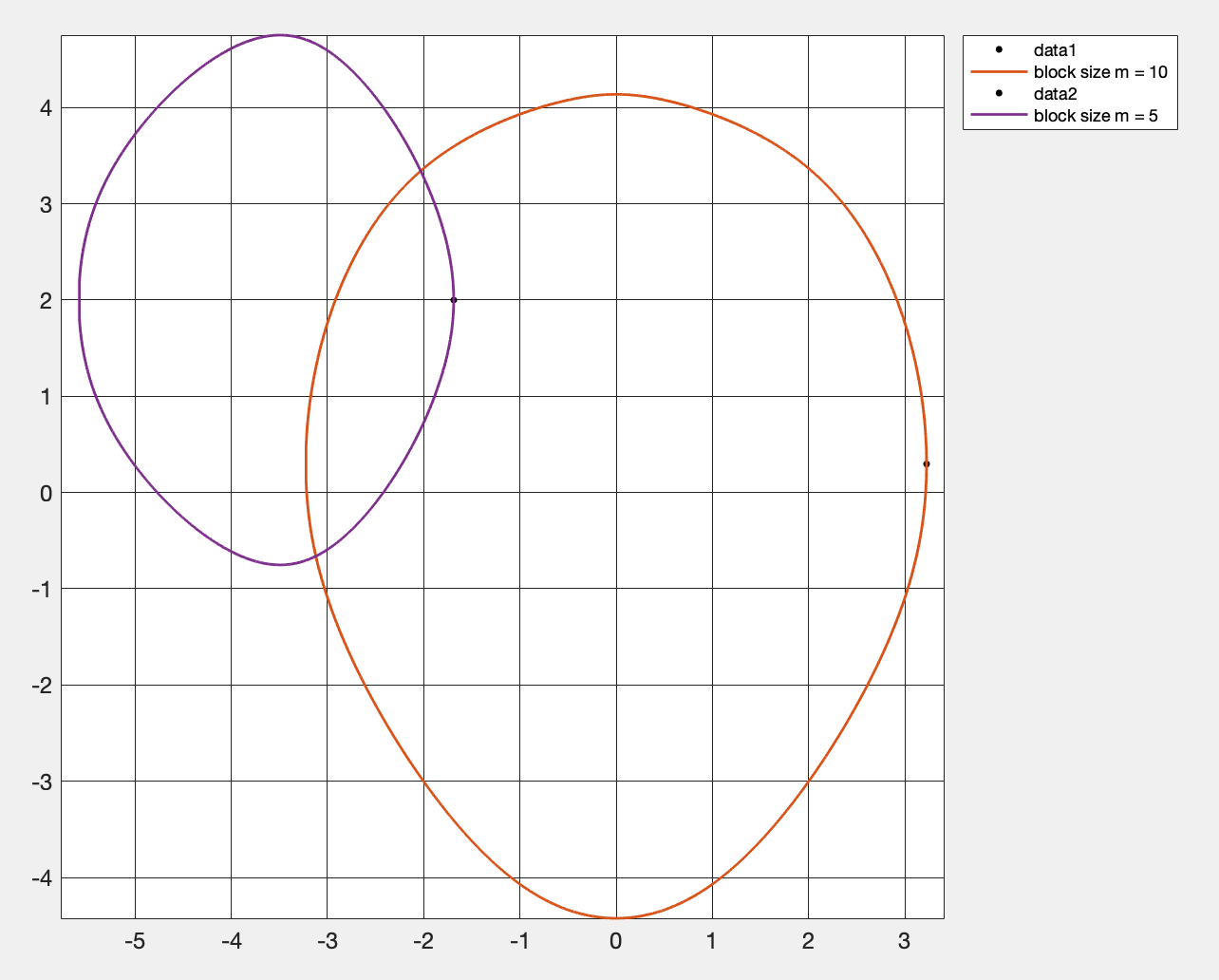} \\[-1mm]
Figure  2: Path following FoV output of our dense test matrix $A$, plotting the FoV boundary of each diagonal block  for $0 \leq t < 2 \pi$ separately
\end{center}

\vspace*{-6mm}
\begin{center}
\includegraphics[width=95mm]{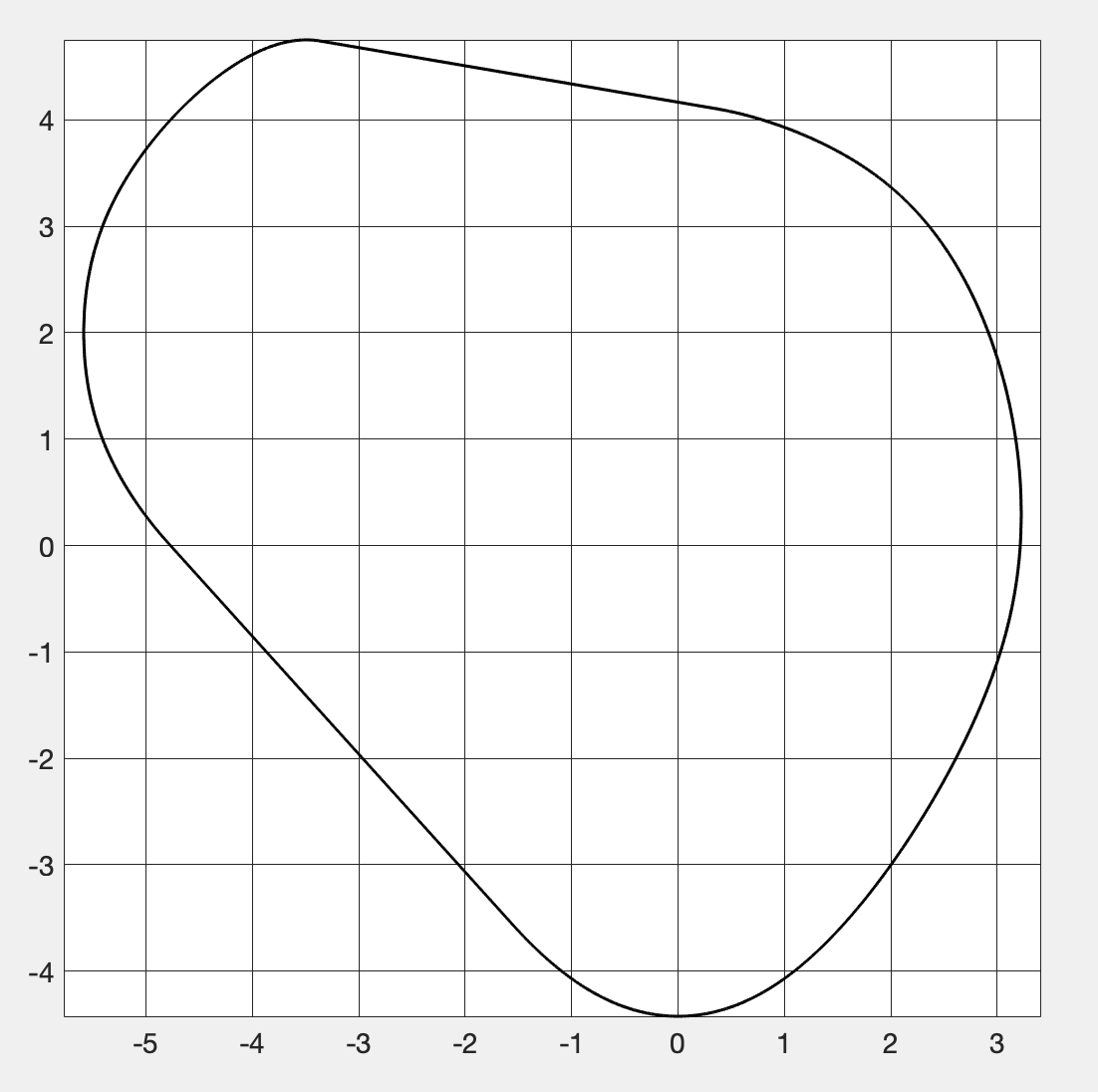} \\[-1mm]
Figure  3: Displaying the FoV boundary curve of $A$ completely, found via a path following method from the individual FoV point data of $A$'s diagonal blocks in Figure 2 and then  extracting  their convex hull
\end{center}

\vspace*{-2mm}
For almost a century mathematicians and physicists have studied the eigenvalue curves of matrix flows $A(t)$. Hermitean matrix flows were studied by Hund \cite{FH1927},  by von Neumann and Wigner \cite{NW29}, by Johnson \cite{J78}, by Dieci et al \cite{DE99,DF01,DP03,DPP13} and others, see \cite[section 5]{FUDecomp}. 
More specifically,  coalescing eigencurves and eigencurve crossing points were studied and computed in \cite{DPP13} using a Newton method based local optimization algorithm.\\[1mm]
Next we display the eigencurves of our chosen block-decomposable test matrix $A$ and then use the block-size assignment algorithm of \cite{FUCoalesc} to find $A$'s hidden block-structure to locate the lead changes in the extreme eigencurves of ${\cal F}_A(t)$.\\[-6mm]
\begin{center}
\includegraphics[width=115mm]{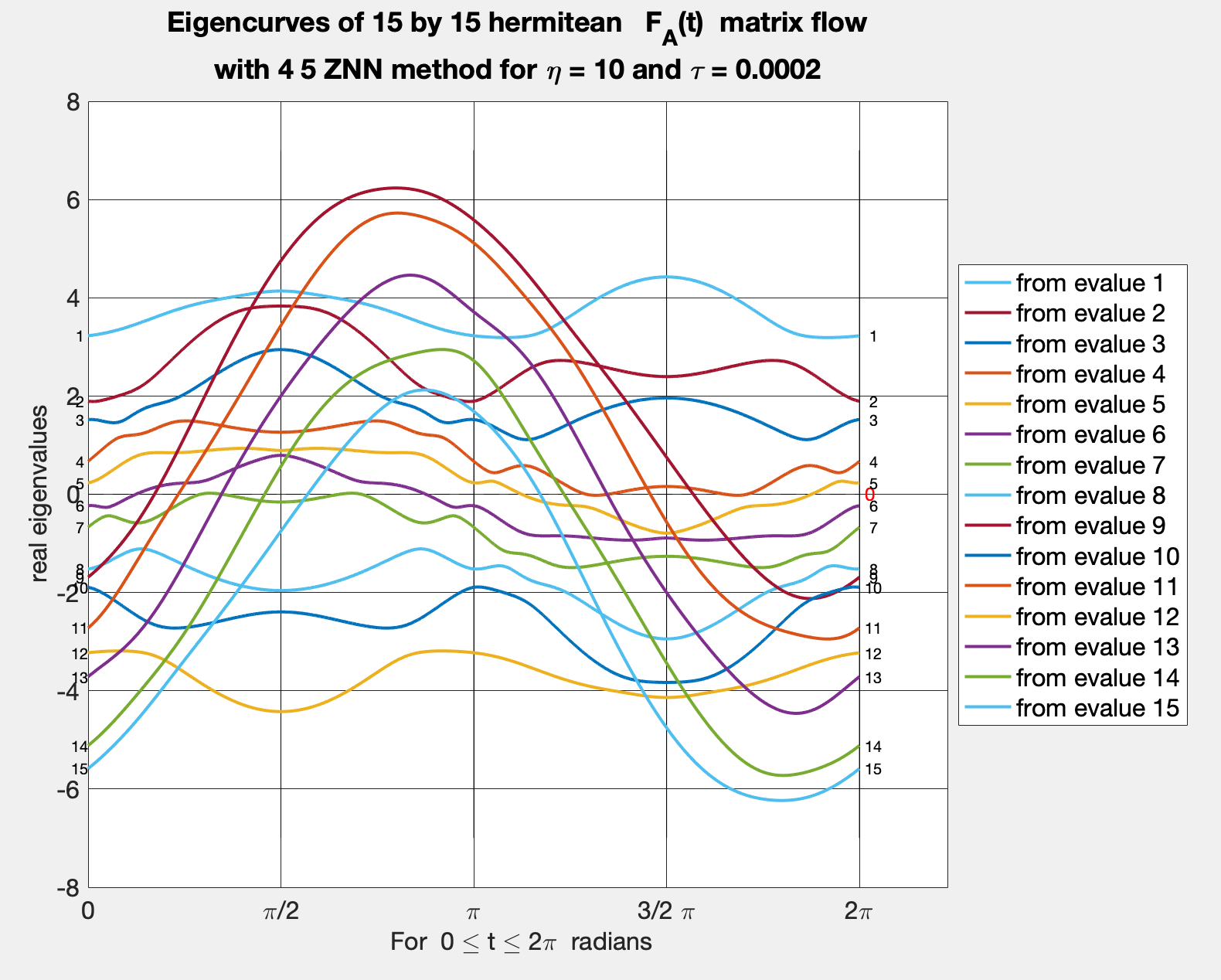} \\[-1mm]
Figure  4: Eigencurves for $A$ and $0 \leq t \leq 2\pi$
 \end{center}
 
\vspace*{-1.5mm}
In Figure 4, note that one eigencurve group  crosses other eigencurves freely, while another eigencurve group veers off hyperbolically from some of the eigencurves. This coalescing and avoidance behavior was first discovered and interpreted by von Neuman  and Wigner  in 1929 for hermitean matrix flows, see \cite{NW29}. Eigencurves that veer off from each other and do not cross are associated with the same diagonal block for a unitarily similar block representation of the flow and those that cross are associated with different diagonal blocks of a  unitarily similar block representation. The simple geometric coalescing eigencurve algorithm of \cite{FUCoalesc} shows that the visually dense matrix $A$ is unitarily decomposable into two diagonal blocks of dimensions 10 and 5 as shown  in Figure 5 where 2 colors suffice to depict the respective eigencurve groups of our 15 by 15 dense test matrix $A$.\\[-4.5mm] 
\begin{center}
\includegraphics[width=95mm]{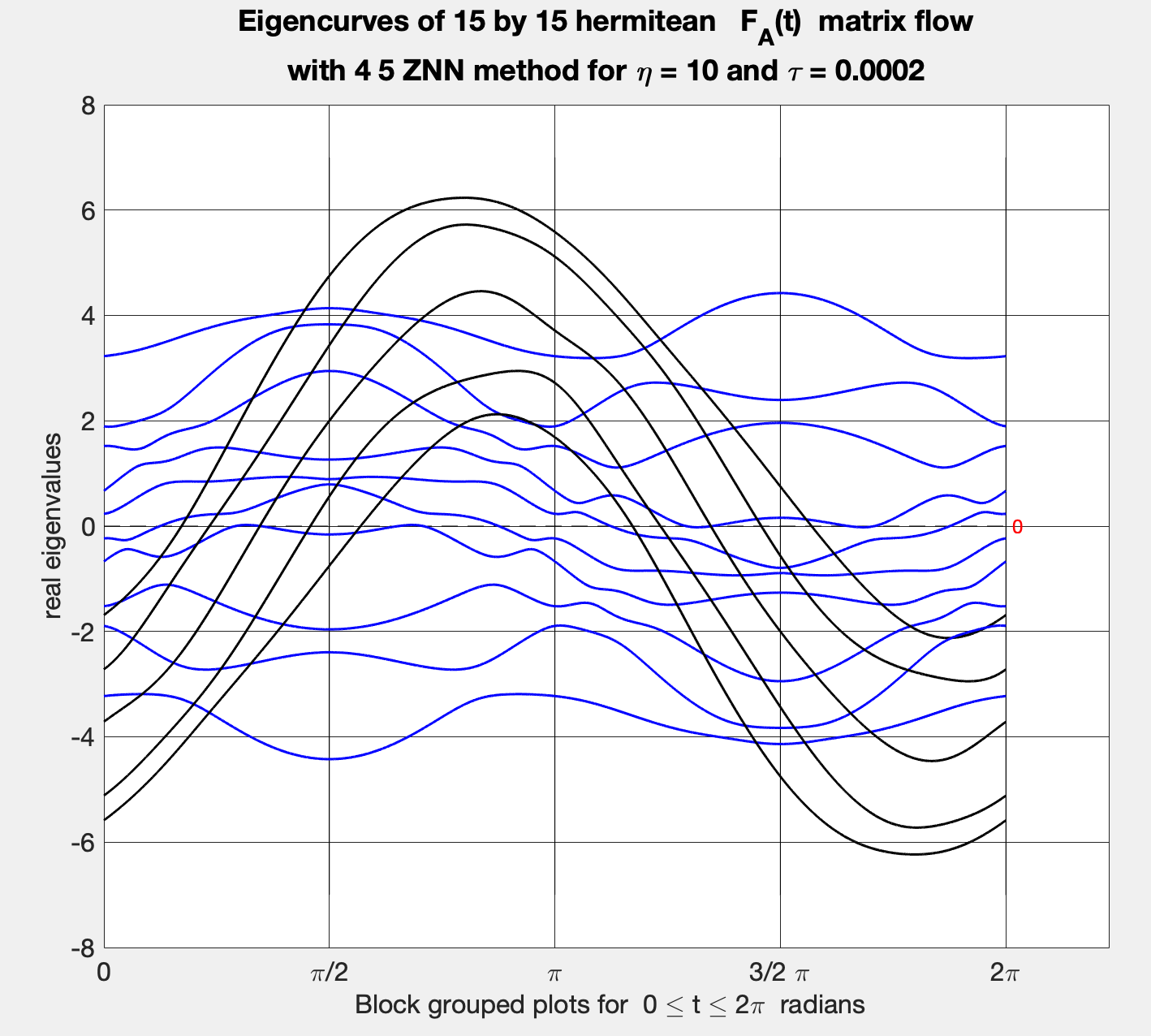} \\[-1.5mm]
Figure  5: {Separated block} eigencurves {(in black and blue)} for $A$ and $0 \leq t \leq 2\pi$
 \end{center}
  It is easy to assess the 'lead changes' among the 'maximal' eigenvalue curves of ${\cal F}_A(t)$ for our test matrix $A$ in Figures 4 and 5.  The  first lead change  occurs just a little before $\pi/2$ {(blue to black)} or at around $ 80 ^o$ when rotating  the Bendixson rectangle or the positive real axis in Figure 3 counterclockwise around the origin. The  lead among the maximal eigenvalue curves reverts back in Figure 5 to the previous eigencurve leader {(black to blue)} at an angle nearly halfway between  $\pi$ and $3/2 \ \pi$ or at around $ 180 ^o + 40 ^o = 220 ^o$ degrees. Both of these  values corroborate well with the two straight line $\partial F(A)$ parts of Figure 3 if we realize that the eigencurve crossing  angles above match the angles of the perpendiculars to the two straight line segments of the FoV boundary curve for $A$ in its upper region and  lower left region  in Figure 3.\\[-6mm]

 \section{The ZNN Path Following Method for General Indecomposable and Decomposable Matrix Fields of Values}
 
 \vspace*{-2mm}
 \subsection{Complexity gains for path following methods}

We first assess the theoretical CPU time gains from using  path following 'divide and conquer' matrix algorithms of complexity $O(n^3)$ when applied to properly decomposable $n$ by $n$ matrices $A$ and matrix flows $A(t)$.\\
Once we  know of a unitary block-diagonalization into smaller blocks of sizes $m_i < n$ with $\sum_{i} m_i = n$  for an  $n$ by $n$ matrix $A$ or a matrix flow $A(t)$,  we can speed up  any $O(n^3)$ unitarily invariant computing process on $A$ or $A(t)$ to the sum of several $O(m_i^3)$  processes with $m_i < n$ for the diagonal block dimensions $m_i$ of a unitarily similar block-diagonal  representation of $A$ or $A(t)$.\\  
If $A$ or $A(t)$ cannot be unitarily block reduced, then  such numerical computations  will still require the full dense matrix $O(n^3)$ effort.\\ 
At the other extreme, if $A$ can be diagonalized by unitary similarity into 1 by 1 diagonal block form, i.e., if $A$ is normal, then any  unitarily invariant  $O(n^3)$ process  on $A$ takes only $O(n)$ effort after -- for example -- an $O(n^3)$ QR based eigenvalue diagonalization of $A$ was obtained, to compute the desired result.  What happens in between these extremes? \\[1mm]
For decomposable $n$ by $n$ matrices and any $O(n^3)$ numerical process we look at the maximal block size $1 \leq m \leq n$ of a given matrix $A$ once its block sizes have been computed via {\tt  decHKflowFoV.m} in \cite{FUDecompMatrFoV}.\\ 
If $A$ has 1 by 1 diagonal blocks except for its one $m$-dimensional block, then the full $O(n^3)$ process becomes one of smaller complexity $O(m^3) = (m/n)^3 O(n^3)$ or $O(m^3) = \alpha^3 O(n^3)$ for $\alpha = m/n$ with $0< \alpha < 1$ as  the remaining one-dimensional blocks require almost no work. This is the best possible scenario when the largest block-diagonal dimension that can be unitarily achieved for $A_{n,n}$  is $m < n$. \\
The worst operations savings scenario happens when the block-diagonalization of $A$ has {\tt floor}($n/m$) blocks of the same maximal size $m$ and  only one additional smaller block, needed  so that the individual block dimensions $m_i$ add up to $n$. Until $m = n/2$ there can be at most {\tt floor}($n/m$) maximal size $m$ blocks. Thus for $m \leq n/2$ the total operations cost of dealing with {\tt floor}($n/m$) maximal size $m$ diagonal blocks is around\\[-2mm] 
 $$\text{{\tt floor}}(n/m) O(m^3) \leq n/m \cdot O(m^3) = n/m \cdot (m/n)^3 O(n^3) = (m/n)^2 O(n^3) = \al^2 O(n^3)$$
 
 \vspace*{-0mm}
  operations for $\al = m/n \leq 1/2$ or $m \leq n/2$. If $m > n/2$ is the size of the largest diagonal block for $A$, then in the 'worst case' there would be one additional indecomposable block of size $n-m < m$ resulting in $O(m^3) + O((n-m)^3)$ necessary operations. Below is a graph of the extreme bounds for the operations cost for a matrix that decomposes unitarily into diagonal blocks of maximal size $m = \alpha \cdot n$ with $0 < \alpha \leq 1$. This graph describes the CPU time gain situation in terms of $0 < \alpha = m/n \leq 1$.\\[-5mm]
\begin{center}
\includegraphics[width=108mm]{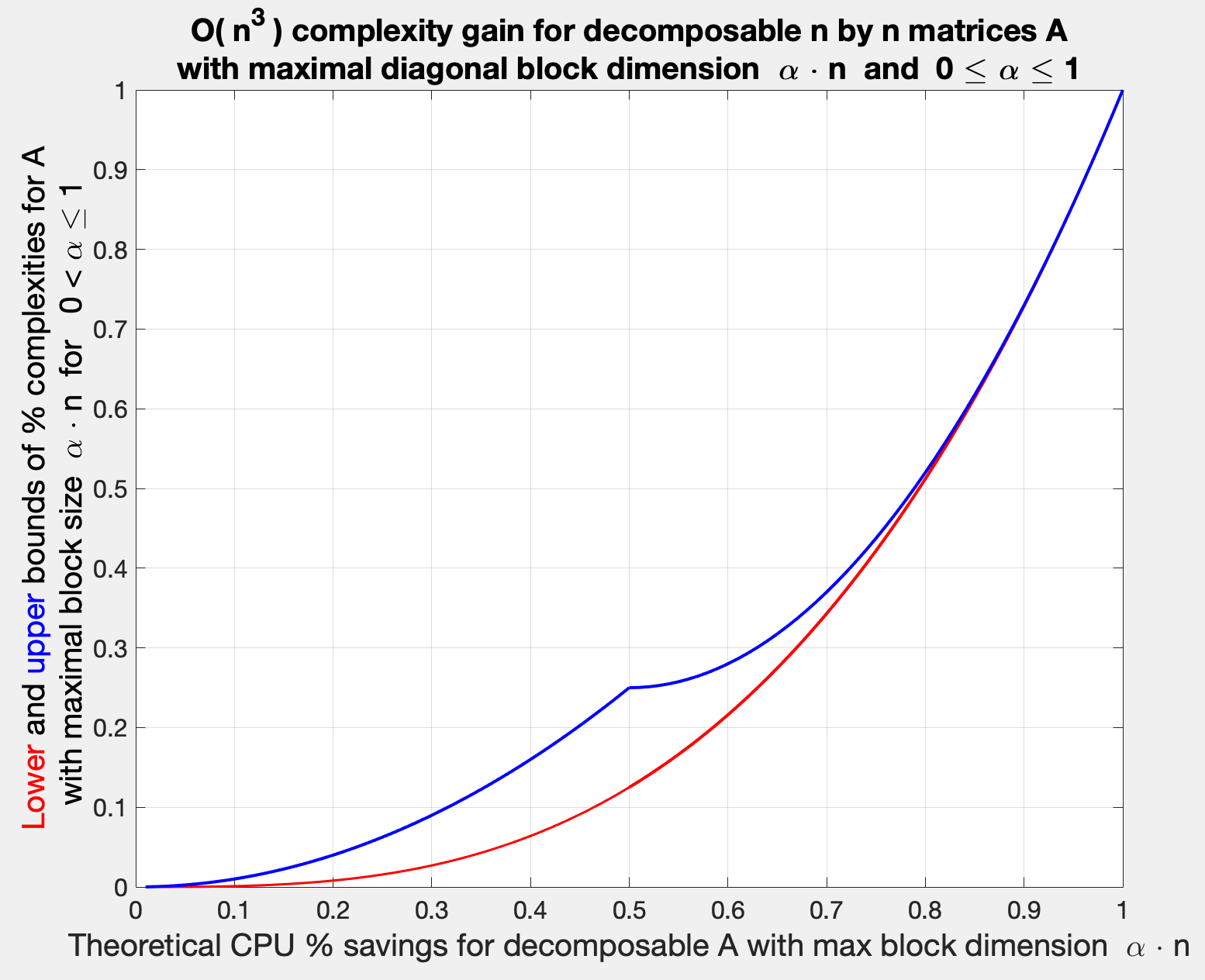} \\[-1.5mm] 
Figure  6 \ \ \ Complexity gains
\end{center}

\vspace*{-2mm}
The largest gap in possible FoV operations counts for decomposable matrices with maximal block dimension $m$ and unitarily invariant computations  occurs at $\alpha = m/n = 0.5$ or when $m = n/2$, see  Figure 6. The theoretical CPU time savings at $ m = n/2$  lie between  87.5 \%  and 75 \% compared to the full $O(n^3)$ effort. Sizable runtime savings of at least 50\% occur  for all $m < 0.7 n$ as can be read off  Figure 6.\\
 Recall the 15 by 15 block-diagonalizable test matrix $A$ of Section 2 whose largest diagonal block had size 10 by 10.There $\al = m/n = 10/15 = 0.667$. Figure 6 indicates that for any $n$ and $\al = 0.667$ the CPU time savings of block-wise computations instead of full $n$ by $n$ matrix computations would be around 67 \%. 

\subsection{Comparing of  path following methods with the global Francis QR eigen algorithm for finding the FoV of indecomposable  and decomposable matrices}

How frequent are unitarily decomposable matrices?  15 non-hermitean matrices in  Matlab's test matrix gallery have difficult to compute eigenstructures and eigenvalues.\\[0.5mm]
 The {\tt 'hanowa'} matrix in this gallery is non-hermitean and normal for all dimensions. Thus it is  unitarily diagonalizable and  benefits from our block-diagonalization algorithm. We note that five out of these 15 matrices or one third are unitarily block-diagonalizable.\\[0.5mm]
  The {\tt 'redheff'} matrix block-diagonalizes for every $n$. Its maximal block size $m$ is usually near its dimension $n$ and there are generally a few additional blocks  of relatively small dimensions such as  2 by 2 or 1 by 1. The FoVs of  Redheffer's small diagonal blocks seem to be  contained within the FoV of its largest diagonal block as depicted below for $n = 27$.\\[-7mm]
\begin{center}
\includegraphics[width=106mm]{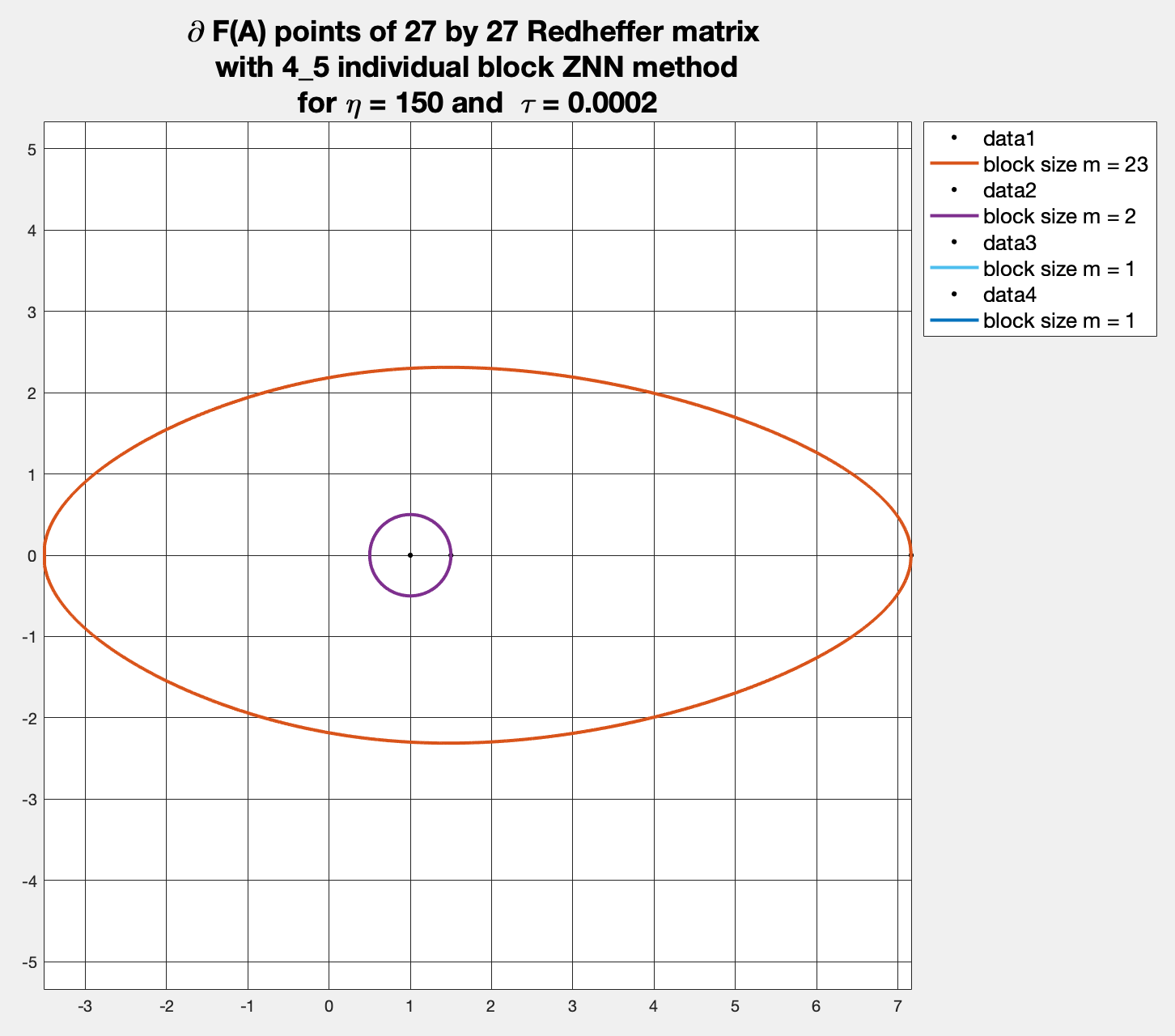} \\[-1mm] Figure  7 \ \ \ Redheffer matrix block FoVs 
\end{center}
Matlab's {\tt 'clement'} matrices of variable sizes $n$ by $n$ always block-diagonalize into two almost equal sized blocks.\\[-10mm]
\begin{center}
\includegraphics[width=104mm]{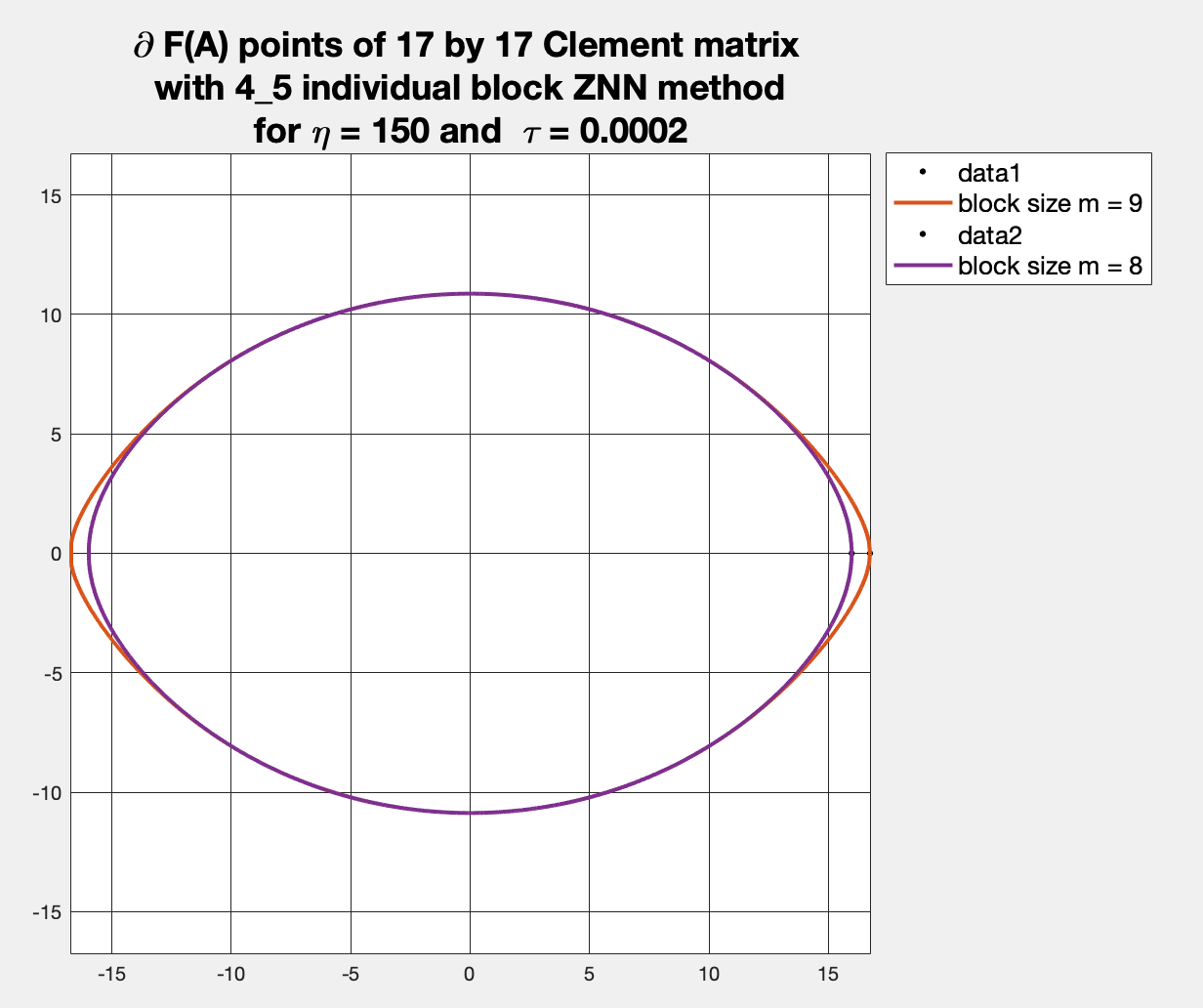} \\[-1mm] Figure 8  \ \ \ Clement matrix block FoVs
\end{center}
And their two diagonal block FoVs  are contained within each other.\\[-0.5mm]
Matlab's $n$ by $n$ gallery matrix {\tt 'binomial'}  unitarily block-diagonalizes into  maximally 4 by 4 blocks when $n$ is even and into maximally 2 by 2 diagonal blocks when $n$ is odd, with additional 1 by 1 blocks so that the block dimensions add to $n$.  \\[-7.3mm]
\begin{center}
\hspace*{4mm}\includegraphics[width=116mm]{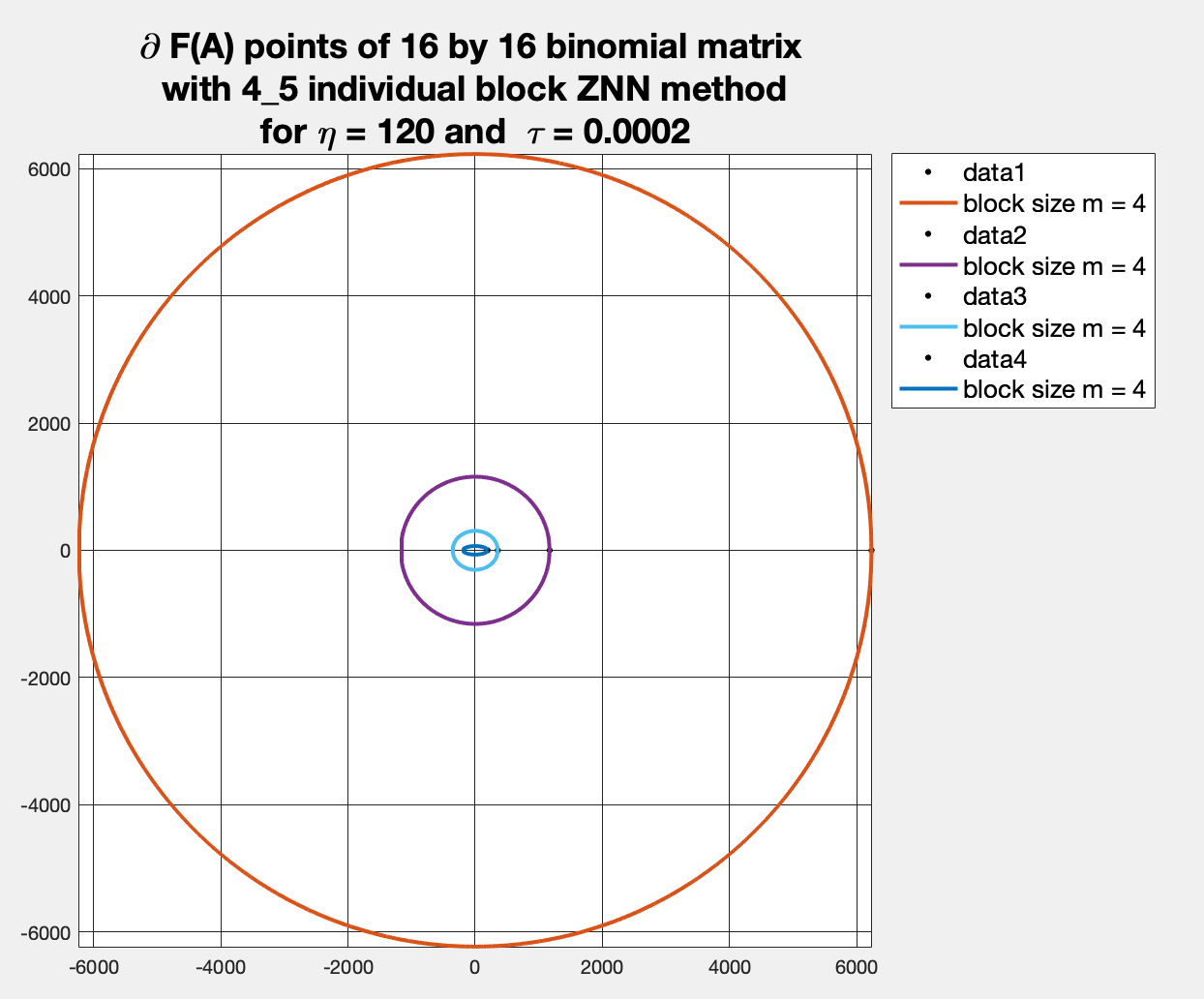} \\[-1mm] 
Figure 9 \ \ \ Binomial matrix separate block FoVs
\end{center}

\vspace*{-1mm}
 Again the  FoVs of the diagonal blocks are nested within each other.\\[2mm]
 The Matlab gallery matrix $A$ = {\tt 'invol'}  is not normal and it achieves uniform maximal size 2 by 2 block-diagonalizations for all dimensions $n$, with one extra 1 by 1 diagonal block for odd $n$. Its separate diagonal block FoV curves seem circular except for the smallest ones. They all seem to be centered at zero and the radii of the larger ones grow rapidly as $n$ increases. When $n = 15$ for example, the maximal radius of its largest 2 by 2 diagonal block FoV curve  measures around $6\times 10^{10}$ units from its center at zero. Its second largest FoV  curve is also nearly circular with center near zero and it has the approximate radius of $3 \times 10^8$.\\[1mm]
 More than a quarter of relevant Matlab gallery  test matrices  are non-normal and three of them or 20 \% can be unitarily  block-diagonalized with maximal block dimensions $m = \al \cdot n$ and $\al < 0.25$. This indicates that their $\partial F(A)$ points can be computed via any path following algorithm and even block QR eigen computations in 10 \% or less of the CPU time of the global $n$ by $n$ Johnson QR method at $O(n^3)$ cost according to the speed gain graph in Figure 6. \\[1mm]
 All of our computations and plots in this paper were obtained using our Matlab m-files for unitary  block-diagonali-\\
 zation in \cite{FUMatrixDecomp} and for the field of values evaluations in  \cite{FUDecompMatrFoV}.\\[2mm] 
Ours is an unusual approach here. Recall that Francis' backward stable QR algorithm and Matlab's {\tt eig} m-file 'diagonalize' every complex matrix  $A$ by finding a non-singular 'eigenvector matrix' $V$ and a diagonal 'eigenvalue matrix' $D$ so that $\tilde A \cdot V = V \cdot D$ with $\tilde A \approx A$, i.e., they solve an adjacent eigenproblem  with $\tilde A$ in a backward stable way for $A$.  Instead we find a unitary $U$ and a proper block-diagonal matrix representation $\hat A$ for any complex matrix $A$ precisely, except for rounding errors, so that $A \cdot U = U \cdot \hat A$ if that is possible. In order to find the FoV of  a decomposable static matrix $A$ efficiently via 'divide and conquer' block methods, we rely on our knowledge of matrix flow decompositions for the related hermitean flow ${\cal F}_A(t)$ of $A$. Thus parameter-varying matrix methods help us  to plot fields of values of static matrices $A$ that block-decompose under unitary similarities. In fact  matrix decompositions for parameter-varying matrix flows now enable us for the first time to use fast eigencurve path following methods such as ZNN without concern for potential eigencurve crossings. Remember that eigencurve crossings do not occur with indecomposable hermitean matrix flows due to Hund and von Neumann and Wigner \cite{FH1927, NW29}. Consequently there are no eigencurve crossings and thus no eigenvalue 'lead changes' within ${\cal F}_{\hat A_j}(t)$ for\linebreak

\vspace*{-3.7mm}
 any indecomposable diagonal block $\hat A_j$ of $\hat A = U^*AU$.\\[1mm]  
The FoV  codes    in \cite{FUDecompMatrFoV} first decide on the unitary decomposability of $A$. If $A$ is unitarily decomposable into  $r$ diagonal blocks $\hat A_j$ for $j = 1,...,r > 1$,  each diagonal block $\hat A_j$  is checked first for normalcy. If normal, then the FoV of $\hat A_j$  is the convex hull of $\hat A_j$'s eigenvalues.  In this case we add the eigenvalues of the normal block $\hat A_j$ to the depository of potential boundary FoV points of $A$. If $\hat A_j$ is not normal, our algorithm computes discrete  FoV boundary points of $\hat A_j$ via ZNN using  a convergent look-ahead finite difference formula of type {\tt k\_s = 4\_5}. It then adds the computed FoV boundary points  of $\hat A_j$ to the depository of potential  FoV boundary points. Upon termination for all  individual blocks, i.e.,  when we know  all sub-FoVs of $A$, our algorithm uses the convex hull algorithm of  Matlab to plot the convex hull of the depository $\partial  F(A)$ point list in $\CC$.  Throughout this paper we use ZNN with a convergent look-ahead   finite difference formula of type {\tt k\_s = 4\_5}  that uses $k+s = 9$ start-up values and has the truncation error order of $O(\tau^{k+2}) = O(\tau^6)$  (see \cite{FUDiscrForm19} for further details) and  the sampling gap $\tau = 0.0002$ while  varying the exponential decay constant $50 \leq \eta \leq 240$ to achieve 15 accurate leading digits for nearly $ 32,000$ discrete FoV boundary points for $A$ that we compute.\\[1mm]
{How can one be sure of the accuracy estimate  '15 accurate leading digits' mentioned 2 lines earlier? Unlike our standard error analysis for static matrix algorithms, ZNN time-varying matrix methods do not lend themselves to explicit error estimates. They cannot use nor need our standard static matrix and vector norm estimations and they are not created to establish any time-varying matrix computational problem or method as backward stable or any time-varying matrix flow $A(t)$ as well or ill conditioned. Continuous and discretized ZNN methods belong to a new uncharted branch of Numerical Matrix Analysis. They are instead governed by the stipulated exponential decay of every entry in the error equation over time. This is due to  ZNN's very design and easily observable. The errors of ZNN matrix methods are, however, still subject to truncation and rounding errors. In discretized ZNN the final accuracy of each of our  finite look-ahead and convergent  difference schemes depends on the sampling gap $\tau$ and its truncation error order, such as $O(\tau^6)$ for formula {\tt 4\_5} above.  This comes from our use of Taylor expansions in the creation of our difference schemes in \cite{FUDiscrForm19} and is born out in practice and  documented in many papers in ZNN's math and engineering literature as well as in \cite[Fig. 1]{FUZNNsurvey} e.g..}\\[2mm]
Below we document several run-time experiments with $n = 250$  and $n = 1000$. For $n = 250$ we limit the maximal diagonal block size $m$ in our test matrices $A_{n,n}$ to $m = 10, 40, 80, 120, 160, 180 <n$ in turn with $m = n =  250$. For $m =10$ and $n = 250$ for example, we  start with  a random entry, block-diagonal matrix $B_{250,250}$ composed of 15 blocks of size $m = 10$, followed by 10 blocks of size $m = 8$ and 5 blocks of size $m = 4$. Then we obscure the block structure of $B$ by forming a  dense test matrix  $A = U^*BU$ which has the identical FoV as $B$ by using a random entry 250 by 250 unitary similarity $U$ on the block-diagonal $B$. Thereafter the actual FoV finding algorithm starts from the dense matrix $A$. First it  retrieves the original block structure of $B$ from $A$ via one hermitean eigenanalysis of  ${\cal F}_A(t_a)$ and a logical 0-1 pattern analysis for a second hermitean flow matrix ${\cal F}_A(t_b) \neq {\cal F}_A(t_a)$  under unitary similarity $V^*{\cal F}_A(t_b)V$ and thereby finds a diagonalizing unitary eigenvector matrix $V$ for ${\cal F}_A(t_a)$ and for $A$ itself. 
This method has been established  in \cite{FUDecomp}. It is  extremely fast,  taking around 0.05 sec of CPU time for general complex matrices $A_{250,250}$ and around 0.9 sec for 1000 by 1000 matrices, both  indecomposable and decomposable ones.  Once we have discovered the  hidden block structure of our test matrices $A$ numerically,  we adapt  our  ZNN eigencode from \cite{FUZNNFoV} to evaluate the eigendata and $\partial F(\hat A_j)$ points of each diagonal block $\hat A_j$ of the block-diagonal realization of $A$ separately. We do this sequentially  for each $j$.
Finally the built-in convex hull  function of Matlab selects the FoV boundary points for $A$ itself.\\[1mm] 
Figure 10 shows  30 individual block FoV curves for  our dense but decomposable example matrix $A_{250.250}$ as specified above. To plot the FoV boundary graph for this $A$, our method uses 31,426 points  out of 942,810 computed potential FoV boundary curve points for $A$'s 30 individual sub-blocks. The whole process takes 20.7  seconds of CPU time.  
Figure 11 shows the block-diagonalization zero-nonzero Matlab  {\tt spy} pattern for ${\cal F}_A(t)$ that is identical to that of the unitarily block-diagonalized version of  $A$ and similar to the block structure  of the original $B_{250,250}$. Using the global Johnson type Matlab {\tt eig} $n$ by $n$ matrix QR based method  takes between 85 and 89 sec for any  maximal diagonal block size $m$ when $n = 250$. The global QR and {\tt eig} based process {\tt wber3FoV2.m} in \cite{FUDecompMatrFoV} computes 32,427 $\partial F(A)$ points and results in a FoV boundary curve whose point coordinates  agree everywhere in their leading 15 digits with those that we have computed via ZNN. Figure 12 finally depicts the computed $\partial F(A)$ curve for $A$ and  our $ m = 10$  example. \\[-6.5mm] 
\begin{center}
\includegraphics[width=135mm]{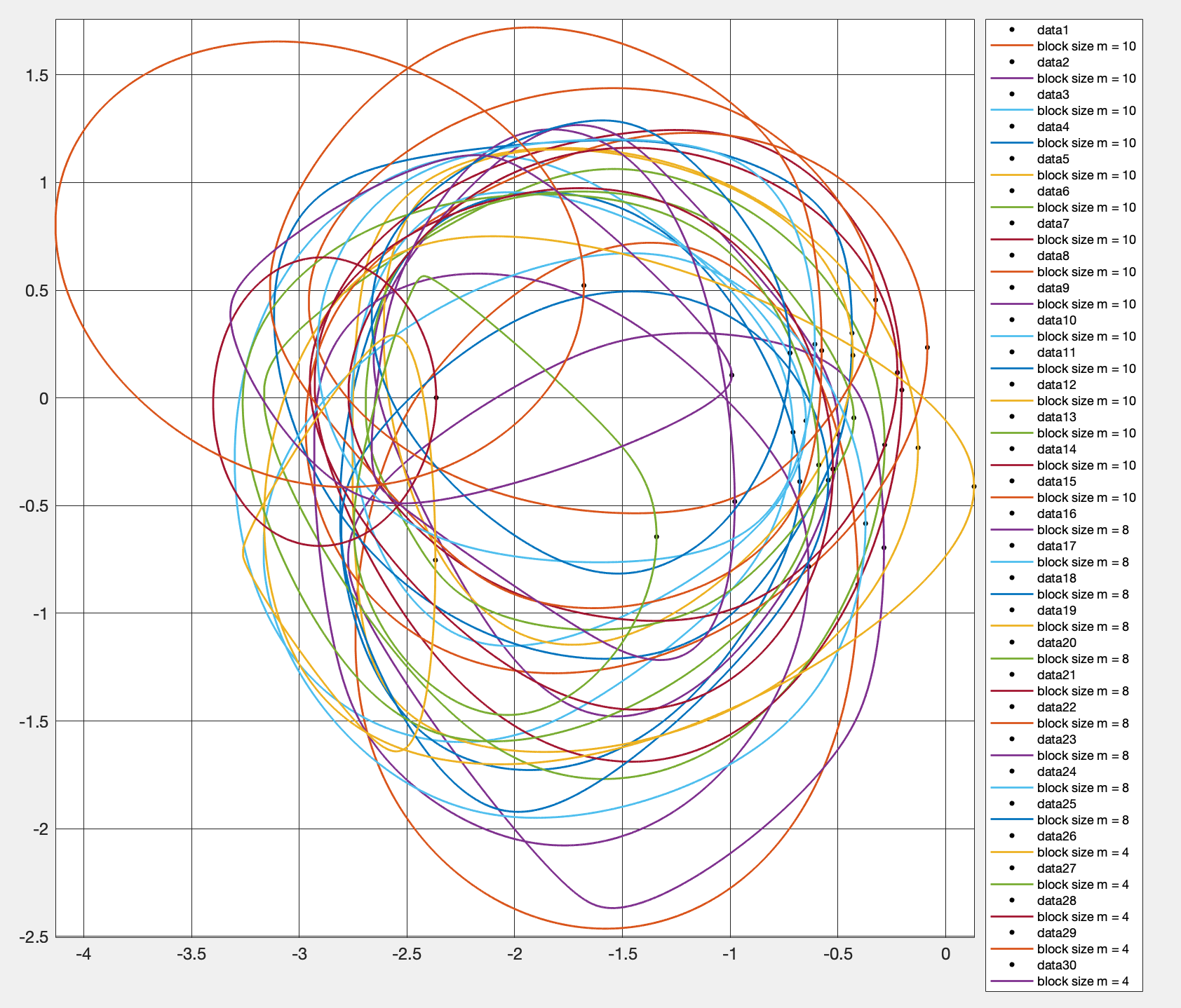} \\[-1mm]
Figure  10:  Showing all 30 diagonal block FoV boundaries of the decomposable  $A_{250,250}$ with their starting points $\cdot$
\end{center}\enlargethispage{40mm}

\vspace*{-8.5mm}
\begin{center}
\includegraphics[width=83mm]{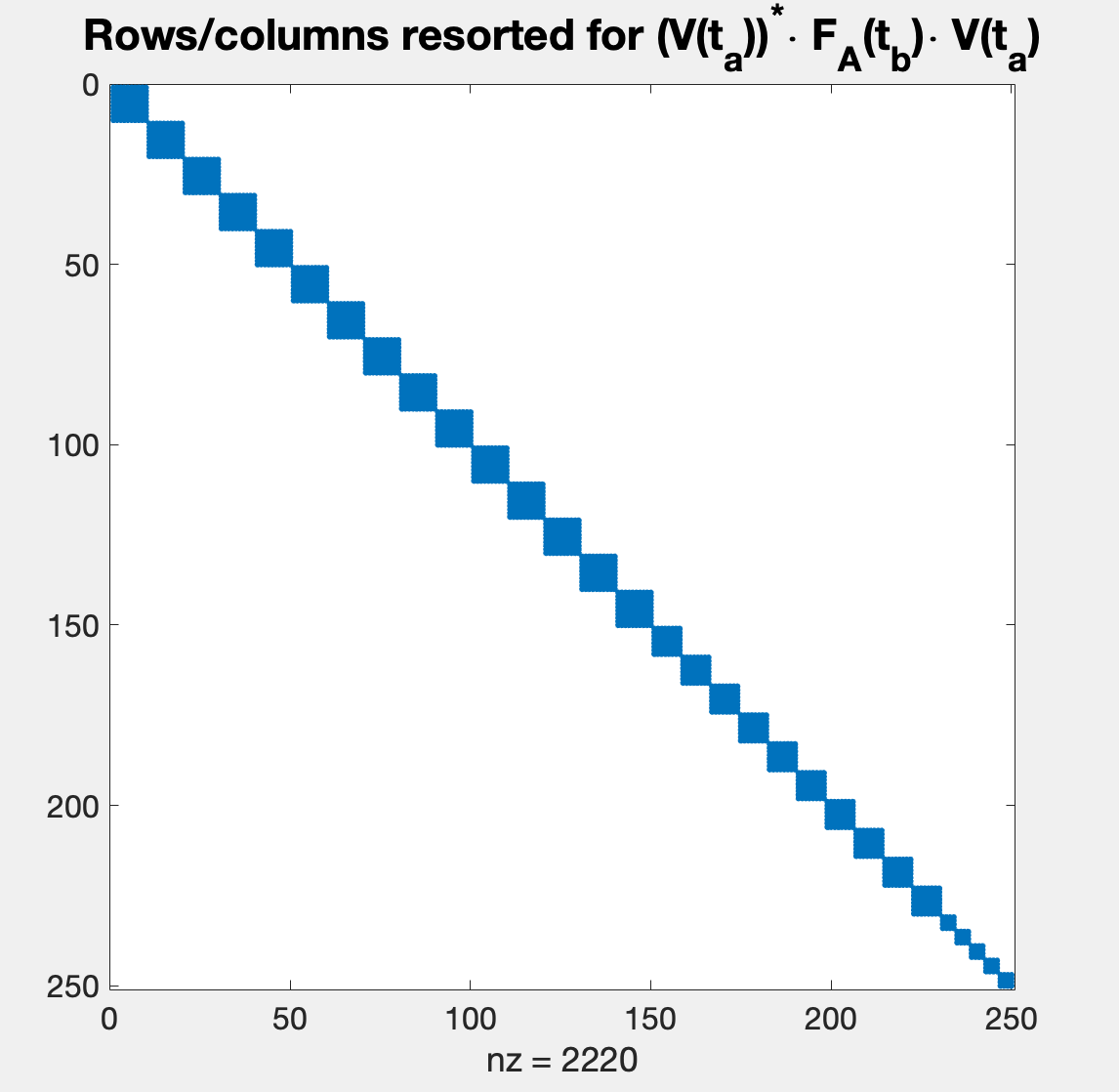} \\[-1mm]
Figure  11: Showing the block-diagonal structure of ${\cal F}_A(t_b)$ under unitary similarity with $V(t_a)$, and thus of $A$\\
\end{center}
\newpage

\begin{center}
\includegraphics[width=90mm]{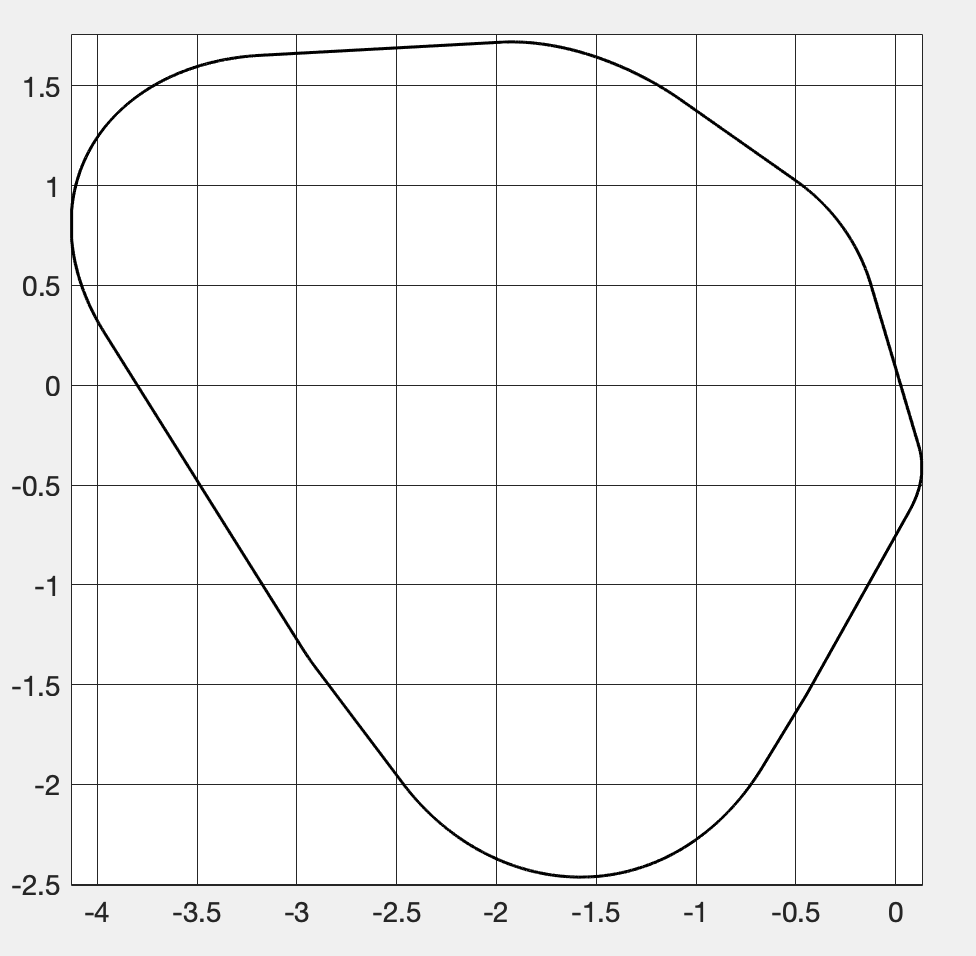} \\[0mm]
Figure  12: FoV boundary curve of the decomposable test matrix  $A_{250,250}$ with maximal block dimension 10,\\ as the convex hull of the partial FoV curve points of Figure 10
\end{center}

The average CPU run times to draw a $\partial F(A)$ curve such as depicted in Figure 12 above  for a given 250 by 250 dense and unitarily decomposable matrix $A$ with increasing maximal diagonal block dimensions $m$   are listed in Table 1. 

\begin{center}\hspace*{-1mm}
\begin{tabular}{c||c|c|c|c|c|c|c}
$\mathbf{A_{n,n}}$ {\bf with}  $\mathbf{n = 250}$&\multicolumn{7}{c} {Max block size $\mathbf m$ for dense unitarily block decomposable  $A_{250,250}$} \\[1mm] \hline \hline
Method  for $\partial F(A)$&10&40&80&120&160&180&250\\ \hline
individual block path following ZNN &20.7  sec & 26  sec& 25  sec&  35 sec&  41 sec&  43 sec &   55 sec\\
 &&\multicolumn{5}{c|}{\emph{ZNN speed-up factors compared to {\tt eig}} :} &\\
 &{\bf 4.3} $\times$&{\bf 3.3} $\times$& {\bf 3.5} $\times$&{\bf 2.5} $\times$&{\bf 2.1} $\times$&{\bf 2} $\times$&{\bf 1.6} $\times$\\ \hline
Johnson QR {\tt eig} based& 89 sec  &86 sec & 88 sec &89 sec &87 sec &85 sec & 87 sec
\end{tabular}\\[2mm]
Table 1 \ \ \ Runtime comparisons; ZNN versus QR  for varying maximal block dimensions $m$ with $n = 250$
\end{center}

\vspace*{-2mm}
The actual run times for our ZNN path following method and FoV computations with fixed $n$ and $m$ may vary by up to around 10 \% up or down, depending on the distribution of block sizes below the maximal block dimension $m$.\\ 
Both our path following ZNN method and the Francis complete  QR  eigendata method can be easily parallelized in Matlab by changing their outermost  'for' loop into a 'parfor' loop that initiates parallel processing. Our sequential 'for' codes  in \cite{FUDecompMatrFoV} end with  {\tt ...for.m} and  their parallelized 'parfor' codes with {\tt ...parfor.m}, respectively. Running both parallelized code versions on the same problem in parallel mode on a 2019 MacBook Pro with a 2.4 GHz 4  core i5 processor and 16 GB of memory  speeds up each of the computations by a factor of around three and the speed advantages of individual block path following methods over global eigendata methods such as Johnson's QR {\tt eig} method stay about the same. The Loisel and Maxwell ODE based path follower in \cite{LM18} can unfortunately not be paralyzed since it does not go through any 'for' loop.\\[1mm]
Our next data table deals with run time comparisons for constructing  the FoV of general square matrices $A_{n,n}$ of sizes $n =250$ and $n = 1000$ that are either dense in the top data rows or that allow a unitary block decomposition into six diagonal blocks of the same proportions for either dimension $n = 250$ or $n = 1000$ in the lower rows.\\[-3mm]

\begin{center}\hspace*{-3mm}
\begin{tabular}{c|r|c||c|r|c|c}
$\mathbf{A_{n,n}}$&\multicolumn{2}{c||}{$\mathbf{n = 250}, \ m = 250$}&$\mathbf{A_{n.n}}$&\multicolumn{3}{c}{$\hspace*{-5mm}\mathbf{n = 1000},\ m = 1000$}\\ 
{\bf dense}&CPU times&\emph{speed-up}&{\bf dense}&CPU time &\emph{speed-up} &$k^x$ ops count\\
&&\emph{vs} {\tt eig}&&&\emph{vs} {\tt eig}& with $k = 4$ 
\\[1mm] \hline \hline\\[-5mm] &&&&&\\[-2mm]
Johnson QR {\tt eig}& 87.2 sec&   &Johnson QR {\tt eig}&3014 sec&    &$k^{2.6}$\hspace*{1.8mm}\\
single ZNN&{ 55.0 sec}&{  1.6}  $\times$&single ZNN&{ 1407 sec}& { 2.1}  $\times$ &${k^{2.34}}$\\[1mm] \hline \hline \hline \\[-4mm]
& \multicolumn{2}{c}{ } & &\multicolumn{3}{c}{}\\[-3mm]
$\mathbf{A_{n,n}}$&\multicolumn{2}{c||}{$\mathbf{n = 250}, \ m = 100$}&$\mathbf{A_{n.n}}$&\multicolumn{3}{c}{\hspace*{-5mm}$\mathbf{n = 1000}, \ m = 400$}\\
{\bf decomposing} into&CPU times&\emph{speed-up}&{\bf decomposing} into&CPU time&\emph{speed-up}& $k^x$ ops count\\ 
{[}100,80,30,30,8,2{]}&&\emph{vs} {\tt eig}&{[}400,320,120,120,32,8{]}&&\emph{vs} {\tt eig}& with $k = 4$ \\[1mm] \hline \hline\\[-5mm]
&&&&&\\[-2mm]
Johnson QR {\tt eig}& 87.2 sec&&Johnson QR {\tt eig}&3014 sec&     &$k^{2.6}$\hspace*{1.8mm}\\
indiv. block ZNN&26.2 sec&3.3 $\times$&indiv. block ZNN&{ 296 sec}& { 10.2}  $\times$ &${k^{1.75}}$\\
\end{tabular}\\[2mm]
Table 2 \ \ \ Runtime comparisons; ZNN versus QR, for dense and decomposing matrices with $n=250$ and $n = 1000$
\end{center}
\vspace*{-3mm}
For real or complex,  dense and unitarily indecomposable matrices $A_{n,n}$ with maximal block dimension $m = n$ and large dimensions $n \gg 25$ our ZNN based path following method is faster than Matlab's Francis implicit multi-shift QR eigendata finder {\tt eig} for plotting FoVs. At $n = 25$ both methods run about equally fast at 1.78 seconds. For $n = 10$ Matlab's fully compiled {\tt eig} needs 0.2 seconds and our ZNN based path finding method computes $A$'s FoV boundary curve in 0.8 seconds. 
Recall from \cite{FUZNNFoV} that ZNN methods reduce matrix eigen computations to solving one hermitean linear systems plus a linear vector recursion at each angle or time step at much lower $O(n^3)$ cost than QR which relies on matrix QR factorizations with their intrinsically  higher $O(n^3)$ costs.\\[0.5mm]
 Progressing from $A_{250,250}$ to $A_{1000,1000}$ with $O(n^3)$ complexity methods for finding FoVs, the data in the right columns  of Table 2 should increase the  operations count  by a factor around the \emph{dimension increase factor $k$ cubed} where $k = 1000/250 = 4$ here, i.e., by around $k^3$ or $4^3$. In fact going to $k=4$ times larger  test matrices  $A$ the operations count increases for the QR based FoV method nominally by less, with an $O(k^{2.6})$ complexity increase in the top right Table 2 entry instead of $O(k^3)$. Note that the operations complexity  increase for single ZNN and FoV plotting is even lower at $O(k^{2.34})$ for indecomposable matrices when going from $n = 250$ to $n = 1000$.  For decomposable matrices $A_{n,n}$ with comparable block decomposition sizes, going from $n=250$ to $n =1000$ improves the ops count even more for our ZNN path finding method from an expected value of $4^3$ to the exceedingly low $4^{1.75}$ ops count as displayed by the bottom right entry of Table 2. \\[2mm]
A dense `artificial' matrix $A_{52,52}$  with complicated  block structure  results in Matlab from the command \\  
{\tt B = blkdiag(-2*eye(2),gallery('forsythe',6), gallery('jordbloc',8,1-1i),...\\
\hspace*{26mm}zeros(3),gallery('hanowa',8),wilkinson(12),hilb(9),...\\
\hspace*{26mm}(gallery('jordbloc',4,1+1i))');} \\
after a dense random entry unitary similarity as $A = Q^* \cdot B \cdot Q$. 
 Matlab's {\tt eig} and our ZNN based path following method (after a block-diagonalization of the dense matrix $A$) plot the $\partial F(A)$  identically as the quadrilateral shown in Figure  13. The respective run times are 2.9 seconds for the {\tt eig} based Johnson QR method and around 2.07 seconds for the individual block  ZNN  method. Note that in the block-diagonalization of the dense test matrix  $A$ only three of its 37 diagonal blocks have dimensions larger than 1 -- one coming from 'forsythe' and two  from the  'jordbloc's of $B$. Otherwise $A$ is  unitarily diagonalizable. The eigenvalues of 'hanowa' lie on the vertical line through $-1 \in \CC$ and the eigenvalues of all other blocks in $A$ are real, denoted by dots on the real axis of Figure 13. The two Jordan blocks of $B$ with the shared eigenvalue $1 - i \in \CC$ generate  concentric circles as sub-fields of values boundary curves for $A$ as is well known. Their radii depend on their dimensions, see Figure 14. \\
 For comparison's sake, the convex hull computation from almost 1 million ZNN computed FoV points for  the decomposing $A_{250,250}$ matrix with $m = 10$ of Table 1 took about 0.08 seconds while the same from around 63,000 FoV points for the unitarily decomposable matrix $A_{52,52}$ takes around 0.05 seconds.

\begin{center}
\includegraphics[width=76mm]{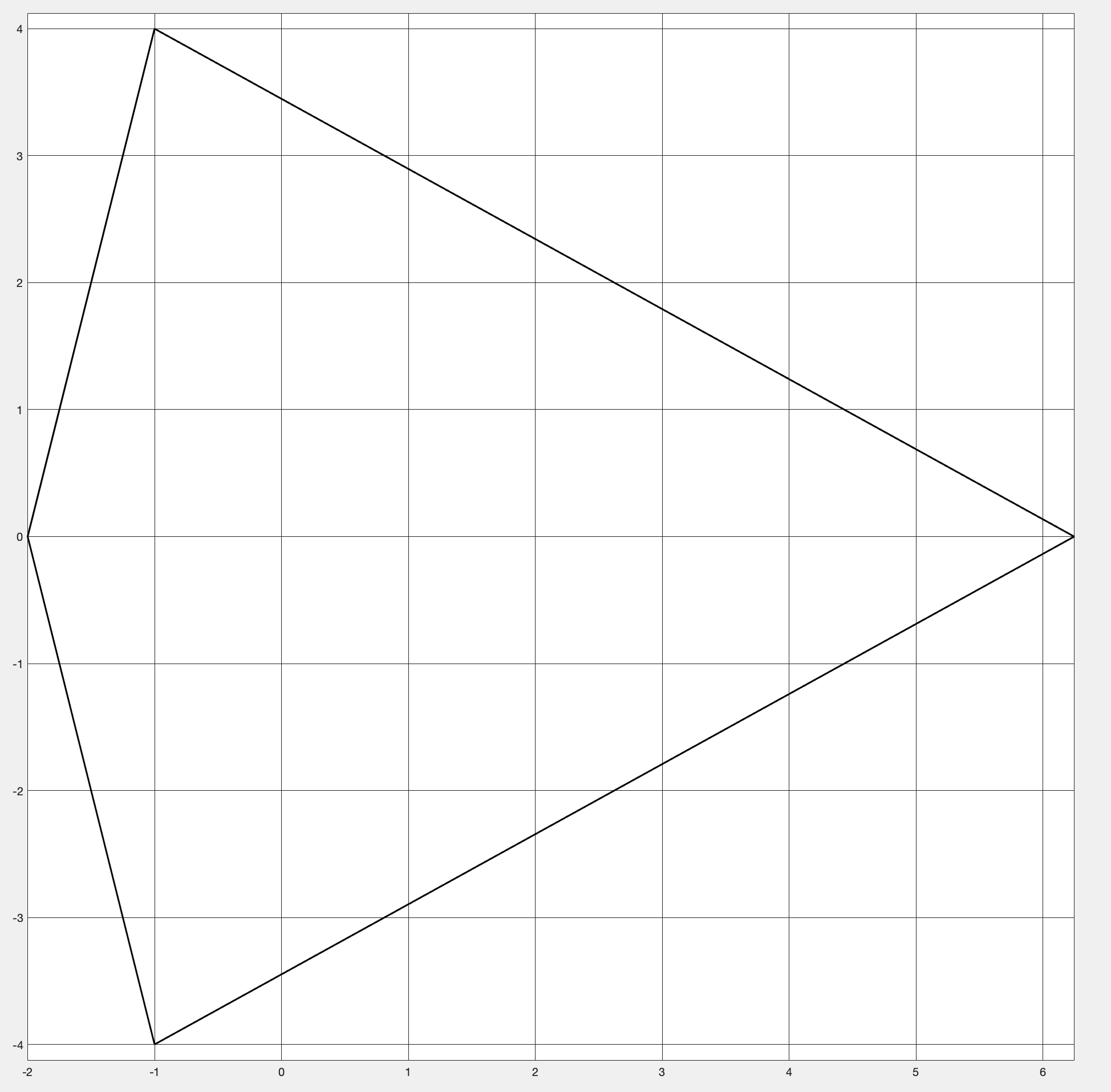} \\[0mm]
Figure 13 : FoV boundary curve of the decomposing test matrix  $A_{52,52}$ 
\end{center}

\vspace*{-3mm}
\begin{center}
\includegraphics[width=88mm]{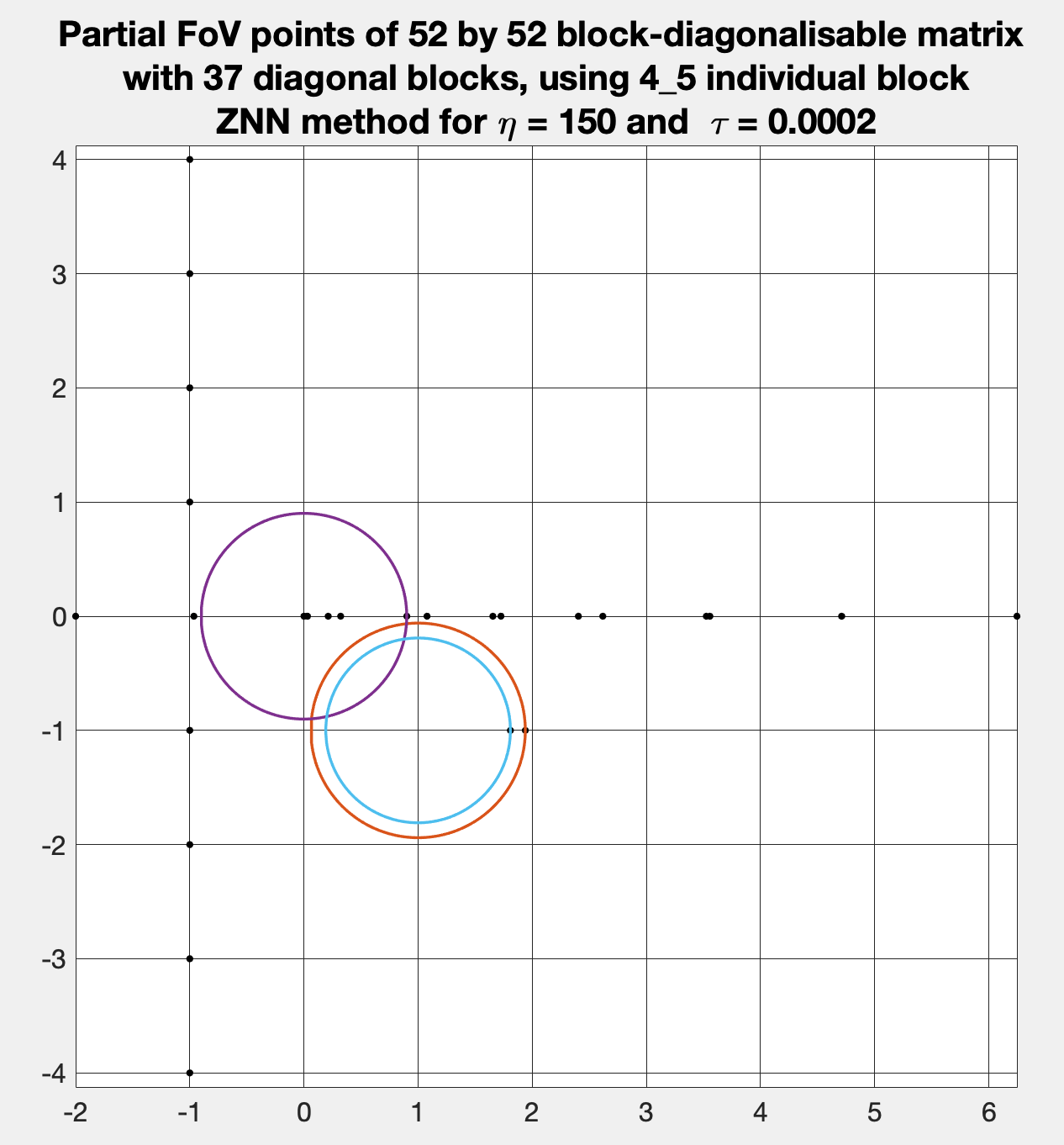} \\[0mm]
Figure 14 : Partial FoV boundary curves of the decomposable test matrix  $A_{52,52}$ depicted under {\tt FOVZNN4\_5aeigshortFoV3for.m} \ in \cite{FUDecompMatrFoV}
\end{center}

A rather debilitating problem arose sometimes unpredictably with computational inaccuracies when replacing $H = H^*$ and $K=K^*$ by the  unitarily similar matrices $V^*HV$ and $V^*KV$ for the diagonalizing matrix $V$ of ${\cal F}_A(t_a)$ that was computed according to the algorithm of \cite{FUDecomp}. Very small un-hermitean  entry errors in $V^*HV$ and $V^*KV$  by quantities of near machine constant \emph{eps} size would occasionally result in Matlab's {\tt eig} function,  interpreting a partial sub-block $H_j$ or $K_j$ of the block-diagonal hermitean matrices $V^*HV$ or $V^*KV$ as non-hermitean and then compute the necessary start-up eigendata  of ${\cal F}_{A_j}(t_k) = \cos(t_k)H_j + \sin(t_k)K_j$ with  nearly correct real parts  and tainted minuscule imaginary parts. Unfortunately however, in this case Matlab's {\tt eig} function's automatic return of ordered  eigenvalues for hermitean matrices  was lost and some computed partial sub FoVs became wild crisscrosses. This made them  useless as  starting values for path following iterations and for finding the final  convex hull FoV for a decomposing matrix $A$. To insure that  all theoretically hermitean diagonal blocks $H_j$ and  $K_j$ are recognized and treated correctly as hermitean in Matlab's {\tt eig} function in the start-up phase for ZNN we  set $\hat H_j = (H_j + H_j^*)/2 = \hat H_j^*$ and $\hat K_j = (K_j + K_j^*)/2 =\hat K_j^*$ before every partial sub-block FoV evaluation and worked with ${\cal \hat F}_{A_j}(t_k) = \cos(t_k)\hat H_j + \sin(t_k)\hat K_j$ instead. \\[2mm]
{Unitarily decomposable matrices appear surprisingly often among the test matrices in the Matlab's matrix gallery. Otherwise  we only know of matrix block diagonalizations in the context of the century old quantum physics problem in Niels Bohr's circle in \cite{FH1927}, \cite{NW29}. There may be many others that so far have been dealt with unawares of their diagonalizable status. Since our test for unitary block-diagonalizability of matrices is so simple, quick and cheap, we advise to check for diagonalizeability periodically when dealing with model matrices in unitarily invariant problems.}

\section{Outlook, backwards and forwards}

Charlie Johnson \cite{J78} may or may not have realized four decades ago that his use of rotating Bendixson  rectangles \cite{Be}  for the field of values problem of a fixed entry static matrix $A$ had  changed the FoV problem fundamentally into a parameter-varying hermitean matrix flow problem for ${\cal F}_A(t)$. The Francis implicit multi-shift QR  method {\tt eig}  is generally considered the most accurate and fastest algorithm for any eigendata computation of square matrices of dimensions not exceeding  $n = 10,000$ by much. Here it evaluates the extreme eigenvalues  of the  hermitean  flow matrices ${\cal F}_A(t)$ of $A_{n,n}$ and records their associated eigenvectors' $x$ actions $x^*Ax \in \CC$   as discrete $ F(A)$ boundary curve points for $A$. Time- or parameter-varying matrix flows $A(t)$ were not part of our Linear Algebra canon nor standard  knowledge then.\\[1mm]
 Even today,  time-varying matrix flow problems have almost exclusively been studied and used in  real-time engineering projects and their computational treatments are just beginning to appear in preprints by numerical analysts and matrix theoreticians. Currently there is only a very thin knowledge base for these problems in  western numerical analysis circles. Yet there are well over 400 papers in engineering journals on Zhang Neural Network (ZNN)  methods for time-varying matrix problems that deal with a wide range  of   on-chip and real-time applications in industry. There is a handful of books on this new, yet mostly theoretically unstudied matrix computational subject, see \cite{FUZNNsurvey} e.g..\\[2mm]
 {The field of values algorithm of Braconnier and Higham \cite{BH96} is a variation of Johnson's QR and Bendixson's rectangle based algorithm that is deliberately less precise to achieve accuracy only 'to ''visual precision" ', see \cite[p. 430, l. -20, -13]{BH96} and it is faster than Johnson's. Instead of QR extreme eigendata evaluations for ${\cal F}_A(t)$ it uses Lanczos-Arnoldi and Lanczos-Chebyshev iterations and is preferrable for sparse matrices. In fact some of their ideas have been incorporated into Matlab's {\tt eigs} m-function \cite[p. 1745]{LM18}.  Braconnier-Higham starts from 16 (or 64) uniformly spaced QR based nearly exact field of values border points for $0 \leq t < 2\pi$ and then uses adjacent FoV boundary point data and Lanczos-plus to fill in the $\partial F(A)$ information of points in between. This method is called 'continuation' and is not path following. The Braconnier-Higham method, however,  is global and handles  decomposing matrices $A$ without any problem. \cite{BH96} gives several references to similar continuation methods for singular value computations \cite{LSH97}, spectral radius estimates and so forth. There is no mention of any other known global FoV algorithm in \cite{BH96} or \cite{LM18}.\\[2mm]
 The first local, i.e., non-global $F(A)$ boundary path following method is apparently in \cite{LM18}. Loisel and Maxwell \cite[sect. 5, 6.2, 7.1, 7.2, 8.3; p. 1733 - 1743]{LM18}  search for ways to solve the FoV problem using $\partial F(A)$ boundary path followers for decomposing matrices and find certain seemingly natural limitations of IVP ODE path following FoV solvers when $A$ is non-normal and derogatory with repeated eigenvalues, i.e., when ${\cal F}_A(t)$ has crossing eigencurves or, equivalently, when $A$ can be unitarily block-diagonalized.} 
Their observations in \cite[last sentence on p. 1743]{LM18} on this issue end in  \\[0.5mm]
\newpage
 
"\emph {... We were also not able to completely analyze nonnormal matrices of type 3 (when $\lambda_{\text{max}}(t)$ can be nonsimple). These limitations to our analysis are to be expected: even for the problem of computing eigenvalues, state-of-the-art eigenvalue solvers fail for some matrices.}"\\[1mm]
  Yes, our current state-of-the art  static matrix eigensolvers such as Francis' multi-shift implicit QR generally  fail for some  matrices $A$ due to seemingly  inherent limitations for matrices with  non-trivial Jordan structures or repeated eigenvalues. Likewise  naive path-following FoV methods must fail  for decomposable  matrices $A$ because potential eigencurve crossings of the associated hermitean matrix flow  ${\cal F}_A(t)$ are non-trivial and not easy to locate.\\[1mm]
   But such 'limitations'  need not be 'expected' for the FoV problem at all as we were told; they are not germane for these problems  -- if we can  enlarge our horizon and understandings. Simple parameter-varying matrix flow ideas have been developed while   studying time-varying matrix problems, see \cite{FUCoalesc} and \cite{FUDecomp} for example. These ideas have now  helped  us to resolve the well 'documented', but in fact  needless FoV  path finder method limitations in a radical new way. We  have shown  how to compute general static matrix FoVs of every square matrix $A \in \CC_{n,n} \text{ or } \RR_{n,n} $, decomposable or not,  -- and how to do so accurately and efficiently -- with  the new and fast ZNN  path-following method when it is  combined with the ultra-fast matrix block-decomposition method of \cite{FUDecomp}. Besides, the matrix and matrix flow unitary decomposition results of \cite{FUDecomp}  that helped us do so here can easily be adapted as  'preconditioner' for other path finding methods and other unitarily invariant matrix problems such as  SVD computations, to least squares problems, to finding the numerical radius or Crawford number of a matrix,  and so forth.\\[1mm]
  Maybe some of the many currently open questions regarding  time-varying matrix flows $A(t)$, see e.g.  \cite{FUDiscrForm19}, \cite{FUZNNFoV}, \cite{FUZNNsurvey}, can be solved in turn by using  deep or simple insights and understandings from  classic fixed entry static matrix analysis, from its concepts, and numerical methods. \\[1mm] 
  I do hope so.\\[-1mm]
  \hspace*{38mm} \underline{\hspace*{60mm}} \\[5mm]
   I gratefully acknowledge the advice and helpful  suggestions from the editor and the referees.

\vspace*{3mm}

\noindent
\centerline{{[} .. /box/local/latex/DecFoVpaper21/DecompMatrFoV3.tex] \quad \today }

\vspace*{5mm}

\noindent
 14 image files :\\[0mm]

F1Blk105folly.png\\
F2Blk105two.png\\
F3Blk105chull.png\\
F4Blk105ecurva.png\\
F5Blk105ecurvblka.png\\
F6CPUtimegaina.png\\
F7Redheff27.png\\[-32.2mm]

\hspace{70mm}\begin{minipage}{70mm}{
F8Clement17.png\\
F9binom16.png\\
F10allblkFoVs250\_30.png\\
F11diag01patt250\_30a.png\\
F12FoV250\_30all.png\\
F13nonsenseA52FoV.png\\
F14B-d52b.png}
\end{minipage}

\end{document}